\documentclass[aap,seceqn,MSNbibl,citesort,dvips]{arximspdf}

\doi{10.1214/09-AAP634}
\volume{20}
\issue{2}
\pubyear{2010}
\firstpage{696}
\lastpage{721}

\makeatletter

\DeclareMathAlphabet\mathcaligr{OMS}{cmsy}{m}{n}

\newtheorem{theo}{Theorem}[section]
\newtheorem{prop}{Proposition}[section]
\newtheorem{lemm}{Lemma}[section]

  \let\sv@tabnotetext\tabnotetext
  \let\sv@tabnotemark@fmt\tabnotemark@fmt
   \long\def\legend#1{{\let\tabnote@indent\leavevmode\sv@tabnotetext[]{}{#1}}}

\makeatother

\begin{document}
\begin{frontmatter}

\title{Normal approximation for coverage
models over binomial point processes}
\runtitle{Normal approximation for coverage}

\begin{aug}
\author[A]{\fnms{Larry} \snm{Goldstein}\ead[label=e1]{larry@math.usc.edu}} and
\author[B]{\fnms{Mathew D.} \snm{Penrose}\corref{}\thanksref{t1}\ead[label=e2]{m.d.penrose@bath.ac.uk}}
\runauthor{L. Goldstein and M. D. Penrose}
\affiliation{University of Southern California and University of Bath}
\address[A]{Department of Mathematics\\
University of Southern California\\
Los Angeles, California 90089-2532\\
USA\\
\printead{e1}} %adresu isvedimo komanda gale!
\address[B]{Department of Mathematical Sciences\\
University of Bath\\
Bath BA2 7AY\\
United Kingdom\\
\printead{e2}}
\end{aug}

\thankstext{t1}{Supported in part by
the Alexander von Humboldt Foundation through
a Friedrich Wilhelm Bessel Research Award.}

% HISTORY:
\received{\smonth{12} \syear{2008}}
\revised{\smonth{8} \syear{2009}}

% ABSTRACT
%
\begin{abstract}
We give error bounds which demonstrate optimal rates of convergence in
the CLT for the total covered volume and the number of isolated shapes,
for germ-grain models with fixed grain radius over a binomial point
process of $n$ points in a toroidal spatial region of volume $n$. The
proof is based on Stein's method via size-biased couplings.
\end{abstract}

% KEYWORDS
%

\pdfkeywords{60D05, 62E17, 60F05, 05C80,
Stochastic geometry, coverage process, Berry--Esseen theorem,
size biased coupling, Stein's method}

\end{frontmatter}

%s1 ###
\section{Introduction}

Given a collection of $n$ independent
uniformly distributed random points in a $d$-dimensional
cube of volume $n$
(the so-called \textit{binomial point process}),
let $V$ denote the (random) total volume
of the union of interpenetrating balls of fixed radius
$\rho$ centered at these points,
and let $S$ denote the number of balls of radius $\rho/2$
(centered at the same set of points)
which are singletons,
that is, do not overlap any other such ball.
These variables are fundamental topics of interest
in the stochastic geometry of coverage processes
and random geometric graphs \cite{Hallbk,Molchanov,Pe,SKM}.

As $n \to\infty$ with $\rho$ fixed
(the so-called thermodynamic limit), both $V$ and $S$
are known to satisfy a central limit theorem (CLT)
\cite{Moran2,Pe,PY1}.
In the present work we provide
associated Berry--Esseen
type
results; that is, we show under periodic boundary conditions
that the cumulative distribution
functions converge to that of the normal at
the same $O(n^{-1/2})$ rate as for a sum of $n$ independent identically
distributed variables, and provide bounds on
the quality of the normal approximation for finite $n$.

Were we to consider instead a
Poisson-distributed number of points,
that is, a~Poisson point process instead of a binomial one,
both of our variables of interest could be expressed as sums
of locally dependent random variables, and thereby
Berry--Esseen %Berry-Esseen
type bounds could be (and have been) obtained
by known methods \cite{AB,GHH,PR,PYStein}. But with
a nonrandom number of points, the local dependence is lost
and the de-Poissonization arguments in \cite{Pe,PY1} do
not provide error bounds for the de-Poissonized CLTs.
The early work of Moran \cite{Moran1,Moran2} on $V$ was
in response to queries in the
statistical physics literature (including the well-known
paper of Widom and Rowlinson \cite{WR}) which specifically addressed
normal approximation of $V$ for nonrandom $n$, and
in general, it seems worthwhile to study the
de-Poissonized setting since in practice one might well
observe the actual number of points, in which case
the conditional distribution of
any test statistic, based on what is observed,
will be over a binomial rather than a Poisson point process.

The variables $V$ and $S$
are just two of a large class
of variables of interest that can be expressed
as a sum, over the $n$ points, of terms that depend only
on the configuration of nearby points in some sense.
General CLTs have been developed
for such variables
\cite{PY1,PeEJP}
and general Berry--Esseen
type results are available in the Poissonized setting
\cite{GHH,PR,PYStein},
but it remains open to provide a generally applicable
Berry--Esseen
type result for such sums when $n$ is nonrandom
(see however \cite{Chat}, which is discussed further in Section
\ref{secgr}).
However, there seem to be good prospects of adapting
the approach of the present paper (which is new in
the geometrical setting) to
a wider class of geometrical sums.

Our approach to normal approximation is based on
Stein's method via size-biased couplings.
Given a nonnegative random variable $Y$ with
positive finite mean $\mu= \mathbb EY$, we say $Y'$ has the $Y$ size-biased
distribution if $P[Y' \in dy] = (y/\mu) P[Y \in dy]$, or more formally,
if
%
%e1.1 ###
%
\begin{eqnarray}\label{formalsb}
\mathbb E[Yg(Y)]=\mu\mathbb Eg(Y')\qquad
\mbox{for bounded continuous functions $f$.}
\end{eqnarray}
The method of size-biased couplings was introduced
by Baldi, Rinott and Stein \cite{BRS}, who used it to develop bounds of order
$\sigma^{-1/2}$ to the normal approximation to the number of local
maxima $Y$ of a random function on a graph, where
$\sigma^2=\mbox{Var}(Y)$. Goldstein and Rinott \cite{GR} extended the
technique to multivariate normal approximations, and improved the
rate to $\sigma^{-1}$ for the expectation of smooth functions of a
vector $\mathbf{Y}$ recording the number of edges with certain fixed
degrees in a random graph.
In \cite{Gold},
the method is used to
give bounds of order $\sigma^{-1}$
for various functions on graphs and permutations.
%represented as a sum of dependent Bernoulli
%random variables, and it is a result
%from \cite{Gold},

Here we shall use Lemma \ref{LGthm} below, which
%is essentially due to the referee and
improves the constant
in a more general result from
\cite{Gold}.
%stated as Lemma that
Loosely speaking, this result says that given any coupling of $Y$ and
$Y'$ on a common space, an upper bound on the distance between the
distribution of $Y$ and the normal can be found which involves
functions of the joint distribution of $Y,Y'$ in terms of (i) the
uniform distance between $Y$ and $Y'$, that is, the $L^\infty$ norm of
$Y-Y'$, and (ii) the
%standard deviation
variance of $\mathbb E[Y-Y'|Y]$.
%which, in some sense, quantify their
%proximity to each other. Hence, it may be possible to
%find good bounds between the distribution of
%$Y$ and the normal if one can closely couple
%$Y'$ to the given $Y$.

%Lemma
% \ref{condsblem} and \ref{sblem}
% illustrates how a
%variable $Y$ can be size-biased, in the settings of interest to us.
% Typically, couplings in
%which $Y'$ and $Y$ are close, such as those
In Section \ref{couplsec} we show how to find a coupled realization of
$Y'$ that is uniformly close to $Y$, for those $Y$ under consideration
here. To do this we show that here the size-biasing amounts to
conditioning the (binomial) number of points falling in a certain
(randomly located) $\rho$-ball to be nonzero, and can be achieved by
modifying at most a single point location to obtain $Y'$ from $Y$, so
that $\| Y'-Y\|_\infty$ is bounded.
%$Y'$ is uniformly close to $Y$.
This construction may be of independent interest, along with Lemma
\ref{condsblem} (a general result on how to size-bias a conditional
probability) and Lemma \ref{LGthm}.
%In Sections
%we provide estimates for the variance of $\E[Y-Y'|Y]$ and
%in Section
% modifying elements of the underlying state space (in our case,
% the set of possible point locations)
%in some minimal manner to obtain a variable with the size-biased
%distribution.

%s2 ###
\section{Results}
\label{secgr}

Let $d \geq1$ and $n \geq4$ be integers. Suppose $U_1, \ldots, U_n$ are
independent random $d$-vectors, uniformly distributed over the cube
$C_n := [0,n^{1/d})^d$ (we write $U_{i}$ rather than $U_{n,i}$ because
the value of $n$ should be clear from the context). Write $\mathcaligr
U_n$ for the point set $\{U_1,\ldots,U_n\}$. For $x, y$ in the cube
$C_n$, let $D(x,y)$ denote the distance between $x$ and $y$ under the
Euclidean toroidal metric on $C_n$. For $x \in C_n$ and $r >0$ let
$B_r(x)$ denote the ball $\{y \in C_n\dvtx D(x,y) \leq r\}$. Let
$B_{i,r}$ denote the ball $B_r(U_i)$. Given $r$, the collection of
balls $B_{i,r}$ form a coverage process (also known as a germ-grain
model) in $C_n$; see \cite{Hallbk,SKM}. Let $\rho> 0$, and define
%
%e2.2 ###
%e2.1 ###
%
\begin{eqnarray}
\label{Vndef}
V &:=& \operatorname{Volume} \Biggl( \bigcup_{i=1}^n
B_{i,\rho} \Biggr) ;
\\
\label{Yndef}
S &:=& \sum_{i=1}^n \mathbf{1} \bigl\{
%(
\mathcaligr{U}_n
% \setminus\{U_i\} )
\cap
B_{i,\rho}
=
\{U_i\}
% \varnothing
\bigr\}.
\end{eqnarray}
Then $V$ is the total covered volume for the coverage process
with $r = \rho$, while $ S$ is the number of
singletons (isolated balls) in the case $r= \rho/2$, and may also be
viewed as the number of isolated points in the geometric graph
on vertex set $\mathcaligr{U}_n$ with distance parameter $\rho$ \cite{Pe}.

Let $Z$ denote a standard normal random variable.
%and let $\SD(\cdot)$ denote standard deviation.
Given a random variable $X$ with $\operatorname{SD}(X) := \sqrt
{\operatorname{Var}(X)} \in(0,
\infty)$,
let $D_X$ denote
% Let $d_K$ denote
the Kolmogorov distance between
the distribution of $X$ (scaled and centered) and that of $Z$, that is,
\[
%D_X:= \sup_{t \in\R} | P[(X- \E X)/\SD(X) \leq t] - P[Z \leq t] |.
D_X:= \sup_{t \in\mathbb{R}} \biggl| P \biggl[\frac{ X- \mathbb
EX}{\operatorname{SD}(X) }
\leq t \biggr] - P[Z \leq t] \biggr|.
\]
%
%probability distributions, i.e. $d_K (F,G) = \sup_{t \in\R}
%|F(t) -G(t)|$.
Our main results provide bounds in the normal approximation
for $V$ and $S$; if $\rho$ is fixed then as $n \to\infty$,
%
%e2.3 ###
%
\begin{eqnarray}
\label{main02}
%d_K ( \LL((V - \E V)/ \SD(V)), \LL(Z) )
D_V
= O(n^{-1/2}) ;\qquad
%d_K ( \LL( (S - \E S)/\SD(S) ), \LL(Z) )
D_S
= \Theta(n^{-1/2}).
\end{eqnarray}
Recall that $a_n = \Theta(b_n)$ means that $a_n = O(b_n)$
and $b_n = O(a_n)$.
We conjecture that the first bound in
%right hand side of
(\ref{main02}) can
be improved to $\Theta(n^{-1/2})$.

To state our results more precisely, we need further notation.
Set $\pi_d$ to be the volume of the unit ball in $d$ dimensions,
that is, $\pi_d: = \pi^{d/2} / \Gamma(1 +d/2)$,
%(see e.g. eqn (6.50) of \cite{Hua}),
and $\phi: = \pi_d\rho^d$.
We say two unit balls \textit{touch} if their closures intersect, but
their interiors do not.
Let $\kappa_d$ (respectively, $\kappa_d^*$)
denote the maximum number of closed unit balls in
$d$ dimensions that can be packed so they all intersect
(respectively, touch)
a~closed unit ball at the origin, but are
disjoint from each other (respectively, have disjoint interiors).
Then $\kappa^*_d$ is
the so-called \textit{kissing number} in $d$ dimensions, which has been studied
for centuries (see \cite{Conway,Zong}).
It is not hard to see
$\kappa_d^*$ is an upper bound for $\kappa_d$,
and in most dimensions it seems likely that $\kappa_d = \kappa_d^*$,
but $\kappa_2 =5$ whereas $\kappa_2^*=6$.
It is known that
$\kappa_3 =
\kappa^*_3 = 12$.
Set $\kappa_d^+: = 1 + \kappa_d$.

Set $\mu_V := \mathbb E[V]$,
$\mu_S := \mathbb E[S]$,
$\sigma_V := \operatorname{SD}(V)$,
and $\sigma_S := \operatorname{SD}(S)$.
It is straightforward to write down formulae for
$\mu_V$, $\mu_S$, $\sigma_V^2$ and $\sigma_S^2$; see
(\ref{0522a}),
(\ref{varV}) and
(\ref{0521a}).

Our first two main results provide nonasymptotic upper bounds
on the Kolmorogorov distance.
%We use the notation $D_X$ for $d_K(\LL((X-\E X)/\SD(X)), \LL(Z))$.
%
\begin{theo}
\label{thm2}
If $ n > 6^d \phi$, then
\begin{eqnarray*}
D_V
%d_K ( \LL(
\leq
\frac{ \mu_V}{5 \sigma_V^2}
\Biggl(
\sqrt{ \frac{11 \phi^2}{\sigma_V} + \frac{5 \sqrt{\eta_V(n,\rho
)}}{\sqrt{n}}
}
+ \frac{2\phi}{\sqrt{\sigma_V}}
\Biggr)^2
% ( \frac{ 4 \volrho^2}{ \sigma_V} +
% ( \frac{\eta_V(n,\rho)}{n} )^{1/2}
% )
\end{eqnarray*}
with
%
%e2.4 ###
%
\begin{eqnarray}\label{etaVdef}\qquad
\eta_V(n,\rho) &:=&
2 \phi^2 \bigl((3^{d} +1)\phi+1\bigr)^2
\nonumber\\
&&{}\times
\biggl( 1 +
(2^d +1) 6^d \phi+
\biggl( \frac{2 n - 6^d \phi}{
n - 6^d \phi} \biggr)
6^{2d} \phi^2
\biggr)
\\
&&{} + 2\phi^4
\biggl(
3 (4^d + 2^d) \phi+
3(4^d) \phi^{2} \biggl(
\frac{2n- 3(2^d) \phi}
{n- 3(2^d) \phi} \biggr)
+ 4 + \frac{2}{n} \biggr).\nonumber
\end{eqnarray}
\end{theo}
\begin{theo}
\label{thm1}
If
%$\sigma_S^{3} \geq6 (n- \mu_S) (\kappa_d +1 )^2$, and
$n > \max( 3^d, 2^{d+1} +1) \phi$, then
\begin{eqnarray*}
%d_K ( \LL(
D_S
\leq
\frac{n- \mu_S}{ 5 \sigma_S^2}
\Biggl( \sqrt{ \frac{ 11 (\kappa_d^+)^2}{ \sigma_S}
+ \frac{5 \sqrt{\eta_S(n,\rho)}}{\sqrt{n}} }
+ \frac{2 \kappa_d^+}{\sqrt{\sigma_S}}
\Biggr)^{2}
% \frac{2(n- \mu_S)}{\sigma_S^2}
% ( \frac{ 4 (\kappa_d+1)^2}{ \sigma_S}
%+ ( \frac{\eta_S(n,\rho)}{n} )^{1/2}
% )
\end{eqnarray*}
with
%
%e2.5 ###
%
\begin{eqnarray}
\label{etaYdef}\hspace*{30pt}
\eta_S(n,\rho) &:=&
2(1+ 2 \kappa_d)^2 \biggl( 1 + (2^d+1) 3^d \phi+
\biggl( \frac{2n- 3^d \phi}{
n - 3^d \phi} \biggr) 9^{d} \phi^2 \biggr)
\nonumber\\
&&{} + \frac{
%( \kappa_d +1)^2
( \kappa_d^+)^2
}{2} \biggl(
\bigl(2^d + 2(3^d) +3\bigr) \phi\\
&&\hspace*{47pt}{} +
(2^{d+1}+1)\biggl( \frac{2n - (2^{d+1} +1) \phi
}{n - (2^{d+1} +1) \phi}
\biggr)\phi^2
+ \frac{4n -2}{n-1}
%4 + \frac{2}{n-1}
\biggr).
\nonumber
\end{eqnarray}
\end{theo}

By using
% \eq{090520b} instead of \eq{090520a}
%in the proof, the
the inequality $(x+y)^2 \leq2(x^2+ y^2)$, the
bounds in Theorems \ref{thm2} and~\ref{thm1}
can replaced by bounds which are simpler, though
less sharp.\vspace*{1pt}
%
%be bounded,
%respectively, by the simpler expressions
%$
% \frac{2 \mu_V}{\sigma_V^2} ( \frac{ 4 \volrho^2}{ \sigma_V} +
% ( \frac{\eta_V(n,\rho)}{n} )^{1/2}
% )
%$
%and
%$
% \frac{2(n- \mu_S)}{\sigma_S^2}
% ( \frac{ 4 (\kappa_d+1)^2}{ \sigma_S}
%+ ( \frac{\eta_S(n,\rho)}{n} )^{1/2}
% ).
%$
% respectively.
%

The next result confirms that
for large $n$, all of
% one might expect
$\mu_V$, $\sigma_V^2$, $\mu_S$ and $\sigma_S^2$
are $\Theta(n)$, so that (\ref{main02})
% and \eq{main0}
follows from Theorems \ref{thm2} and \ref{thm1}.
To provide details we require further notation.

For $0 \leq r \leq2$, write
$\omega_d(r)$ for the volume of the union of two unit
balls in $\mathbb{R}^d$ with
centers distant $r$ apart (see (\ref{omMoran})
for a formula).
Define the integral
%
%e2.6 ###
%
\begin{equation}\label{Jdef}
J_{r,d}(\rho) := d \pi_d\int_0^r \exp( - \rho^d \omega_d(t) ) t^{d-1}\,dt
\end{equation}
and the functions
%
%e2.8 ###
%e2.7 ###
%
\begin{eqnarray}\qquad
\label{gVdef}
g_V(\rho) & := & \rho^d J_{2,d}(\rho) -
(2^d \phi+ \phi^2)
e^{-2\phi} ;
\\
\label{gWdef}
%and
g_S(\rho) & := &
e^{-\phi} - \bigl(1 + (2^d-2) \phi+ \phi^2\bigr)
e^{-2 \phi} + \rho^d \bigl(J_{2,d}(\rho) -J_{1,d}(\rho)\bigr).
\end{eqnarray}
Also, define
$ \eta_V(\rho) : = \lim_{n \to\infty} \eta_V(n,\rho) $
and
$ \eta_S(\rho)
: = \lim_{n \to\infty} \eta_S(n,\rho)$.
Formulae for
these limits
are immediate from the definitions
(\ref{etaVdef}) and (\ref{etaYdef}).
\begin{theo}
\label{thmlimsup}
If $\rho$ is fixed then as $n \to\infty$,
%
%e2.11 ###
%e2.10 ###
%e2.9 ###
%
\begin{eqnarray}
\label{meanlim}
\lim_{n \to\infty} \bigl( 1-n^{-1} \mu_V(\rho) \bigr) & = &
\lim_{n \to\infty} ( n^{-1} \mu_S(\rho)) =
e^{-\phi};
\\
\label{varVlim}
\lim_{n \to\infty} ( n^{-1} \sigma_V^2 ) & = & g_V(\rho) > 0;
\\
\label{varWlim}
\lim_{n \to\infty} ( n^{-1} \sigma_S^2 ) & = & g_S(\rho) > 0
\end{eqnarray}
and
% (with $\kappa_d^+ := 1+ \kappa_d$)
%
%e2.14 ###
%e2.13 ###
%e2.12 ###
%
\begin{eqnarray}\qquad
\label{main2}
\limsup_{n\to\infty}
(n^{1/2}
%(n^{1/2}
D_V
)
& \leq&
%
% \frac{8 \volrho^2(1-e^{-\volrho}) }{g_V(\rho)^{3/2}}
%+ \frac{2 (1-e^{-\volrho}) \eta_V(\rho)^{1/2} }{g_V(\rho)} ;
%
\frac{1 - e^{-\phi}}{5 g_V(\rho)}
\Biggl(
\sqrt{
% (
\frac{11\phi^2}{g_V(\rho)^{1/2}} + 5
\eta_V^{1/2}
% )^{1/2}
}
+ \frac{2 \phi}{g_V(\rho)^{1/4}
}
\Biggr)^2;
\\
\label{main1}
\limsup_{n\to\infty} (n^{1/2}
D_S
)
& \leq&
%
% \frac{8 ( \kappa_d +1)^2(1-e^{-\volrho}) }{g_S(\rho)^{3/2}}
%+ \frac{2 (1-e^{-\volrho}) \eta_S(\rho)^{1/2} }{g_S(\rho)} ;
%
\frac{1 - e^{-\phi}}{5 g_S(\rho)}
\Biggl(
% (
\sqrt{
\frac{11(\kappa_d^+)^2}{g_S(\rho)^{1/2}} + 5
\eta_S^{1/2}
}
% )^{1/2}
+ \frac{2 \kappa_d^+}{g_S(\rho)^{1/4}}
\Biggr)^2;
\\
\label{SLB}
\liminf_{n\to\infty} (n^{1/2}
D_S
)
& \geq& (8 \pi g_S(\rho))^{-1/2} .
\end{eqnarray}
\end{theo}

Theorems \ref{thm2} and \ref{thm1}
are proved in Sections \ref{proof-2}
and \ref{proof-1}, respectively.
Theorem \ref{thmlimsup} is proved in Section \ref{secvar},
where we also derive numerical values for
%$\delta_W(\rho)$ and $\delta_V(\rho)$,
the asymptotic
upper bounds in Theorem \ref{thmlimsup}, for some particular cases.
%Unfortunately, the large values of these constants may limit
%the practical usefulness
%of our results, except for very
%large samples.

\subsection*{Remarks}
The limiting variances in (\ref{varVlim}),
respectively (\ref{varWlim}),
are consistent with those
given by Moran \cite{Moran1,Moran2},
respectively,
Penrose (\cite{Pe}, Theorem 4.14).
Moran and Penrose do not
%appear to
explicitly rule out the possibility that these limiting variances might be
zero, as we do here.
%an issue we deal with in the proof of \eq{varVlim}.

Clearly (\ref{main2}) and (\ref{main1}) imply central limit theorems whereby
both $(V- \mu_V)/\sigma_V$ and $(S- \mu_S)/\sigma_S$
converge in distribution to the standard normal,
thereby providing an alternative to existing
proofs of these central limit theorems \mbox{\cite{Moran2,PY1,Pe}}.
In the Poissonized setting, nonasymptotic
bounds analogous to those in Theorems \ref{thm2} and \ref{thm1}
are given in \cite{PR}
and imply $O(n^{-1/2})$
bounds analogous to
(\ref{main2}) and (\ref{main1}).
In the de-Poissonized setting
considered here,
Chatterjee \cite{Chat} provides
bounds similar to those in (\ref{main2}) and (\ref{main1}),
which hold for general metric spaces,
but using the Kantorovich--Wasserstein distance, rather than
the Kolmogorov distance considered here, and without
providing any explicit constants.
As stated in~\cite{Chat}, ``obtaining optimal rates for the Kolmogorov
distance requires extra work and new ideas.''

Generalizations of our results
should be possible
in many directions.
These include:

\subsection*{More general germ-grain models}
Replace the balls of fixed radius in the description of
$V$ and $S$ by
(independent identically distributed)
balls of
random radius, or more generally, random shapes.

\subsection*{Random measures}
Consider the random measure associated
with $V$ (the Lebesgue measure on the covered region)
or with $S$ (a sum of Dirac measures at
the isolated points), and look at normal approximation for
the random variable
given by the integral of a test function $f$ on $C_n$ with respect
to that measure.

\subsection*{Euclidean distance}
Suppose in the definition of $V$ and $S$, that
the periodic boundary conditions on $C_n$ are dropped, that is, the toroidal
distance $D$ is replaced by the ordinary Euclidean distance.

\subsection*{Nonuniform points}
Consider a sequence of
independent random points $(\mathbf{X}_n)_{n \geq1}$ with a common
density function
$\nu\dvtx \mathbb{R}^d \to\mathbb{R}$. Placing balls of radius $r_n$
around each point of $\mathcaligr{X}_n := \{ \mathbf{X}_1,\ldots
,\mathbf{X}_n\}$,
for some specified sequence $r_n$ tending to zero,
one may define quantities analogous to $V$ and $S$.
When $r_n \propto n^{-1/d}$ this is a rescaling of our
model but allows for nonuniform $\nu$. Our approach
might also provide information
about other asymptotic regimes.

\subsection*{$k$-nearest neighbors}
Let $k \in\mathbb{N}$ and consider the
number of points $U_i$ whose $k$th nearest neighbor in the point
set $\mathcaligr{U}_n \setminus\{U_i\}$ lies at a distance greater
than $\rho$.
The case $k=1$ reduces to $S$.

These extensions generally seem to be nontrivial, and worthy of further
study.

%s3 ###
\section{Lemmas}

The proof of (\ref{main2}) and (\ref{main1}) is based on
the following result. This result
%is an optimized version of one due to the referee, and
improves the constant which would be obtained by applying
the more general
% which is a special case of
Theorem 1.2 of
\cite{Gold} to the particular case of Kolmogorov distance.
\begin{lemm}
\label{LGthm}
Let $Y \ge0$ be a random variable with mean $\mu$ and variance
$\sigma^2 \in(0,\infty)$, and let $Y^s$ be defined on the same
probability
space, with the $Y$-size biased distribution. If $P[|Y^s-Y| \le B] =1$
for some constant $B >0$,
%$B \le\sigma^{3/2}/\sqrt{6 \mu}$,
then
%
%e3.1 ###
%
\begin{equation}\label{090520a}
%d_K ( \LL( \frac{Y-\mu}{\sigma} ) ,\LL(Z) )
D_Y
\le
\frac{\mu}{5 \sigma^2} \Biggl(\sqrt{\frac{11B^2}{\sigma} + 5 \Delta} + \frac{2 B}{
\sqrt{\sigma}}\Biggr)^2,
%& \leq&
\end{equation}
where
$\Delta: =
\sqrt{\operatorname{Var}(\mathbb E[Y^s-Y|Y])}$.
\end{lemm}
\begin{pf}
% Letting
%$$
%$$
Given $z \in\mathbb{R}$ and $\varepsilon>0$, define real-valued functions
$h_z$ and $h_{z,\varepsilon}$
%from $\R$ to $\R$
by
% function $h \in\mathcaligr{H}$ and $\varepsilon>0$,
%let
%
\[
h_z(x)= \mathbf{1}_{(-\infty,z]}(x),\qquad
h_{z,\varepsilon}(x)=\varepsilon^{-1}\int_0^\varepsilon
h_{z}(x-t)\,dt,\qquad
z \in\mathbb{R}.
\]
%
%We may write the Kolmogorov distance between
%$W=(Y-\mu)/\sigma$ and a standard normal $Z$ as
Then with $W:=(Y-\mu)/\sigma$ and
$Z$
denoting a standard normal,
by definition
%
%e3.2 ###
%
\begin{equation}\label{090528a}
D_Y =\sup\{|\mathbb Eh_z(W) - \mathbb Eh_z( Z)|\dvtx z \in\mathbb{R}\}.
\end{equation}
For $\varepsilon>0$,
set
%
%e3.3 ###
%
\begin{equation}\label{090528b}
D_Y^\varepsilon:=\sup\{ |\mathbb Eh_{z,\varepsilon}(W) - \mathbb
Eh_{z,\varepsilon}(Z)|\dvtx
z \in\mathbb{R}\}.
\end{equation}
Fix $z$ and $\varepsilon$,
%$h \in\mathcaligr{H}$,
and let $f$ be the unique bounded solution of the Stein equation
\[
f'(w)-wf(w)=h_{z,\varepsilon}(w)-\mathbb Eh_{z,\varepsilon}(Z)
\]
for $h_{z,\varepsilon}$; see \cite{CS}.
%Then,
With some abuse of notation, let $W^s=(Y^s-\mu)/\sigma$.
Then
% we have
%
%e3.4 ###
%
\begin{eqnarray}
\label{L-BD:KSepsm}
&&\mathbb E [ h_{z,\varepsilon}(W)-\mathbb Eh_{z,\varepsilon}(Z)
] \nonumber\\
&&\qquad= \mathbb E [ f'(W)-Wf(W) ] \nonumber\\
&&\qquad= \mathbb E \biggl[ f'(W)- \frac{\mu}{\sigma}\bigl(f(W^s)-f(W)\bigr)
\biggr]\\
&&\qquad=\mathbb E \biggl[f'(W)\biggl(1-\frac{\mu}{\sigma}(W^s-W)\biggr)\nonumber\\
&&\qquad\quad\hspace*{10.8pt}{} - \frac{\mu
}{\sigma
}\int_0^{W^s-W}\bigl(f'(W+t)-f'(W)\bigr)\,dt \biggr] \nonumber.
\end{eqnarray}
The following bounds on the solution $f$ can be found in \cite{CS}:
%
%e3.6 ###
%e3.5 ###
%
\begin{equation}
\label{L-BD:ChSh04bd}
|f'(w)| \le 1
\end{equation}
and
\begin{equation}
\label{090615a}
|f'(w+t)-f'(w)|
\le (|w|+1)|t|+ \varepsilon^{-1} \int_{t \wedge0}^{t \vee0}
\mathbf{1}_{[z ,z+\varepsilon]} (w+u) \,du.
\end{equation}
Noting that $\mathbb EY^s=\mathbb EY^2/\mu$ by (\ref{formalsb}) with
$g(y)=y$, we find
that
\[
\frac{\mu}{\sigma}\mathbb E[W^s-W]=\frac{\mu}{\sigma^2} \biggl(
\frac
{\mathbb EY^2}{\mu}-\mu\biggr)=1,
\]
and therefore, taking expectation by conditioning, and then
using
% applying the first bound in
(\ref{L-BD:ChSh04bd}), we have
\[
\biggl\vert\mathbb E \biggl\{ f'(W) \mathbb E \biggl[ 1 - \frac{ \mu}{\sigma}
(W^s-W)\Big|W \biggr]
\biggr\}\biggr\vert
\le\frac{ \mu}{\sigma} \sqrt{\mbox{Var}(\mathbb E
[W^s-W|W])}= \frac{ \mu}{\sigma^2}\Delta.
\]

Now, using
(\ref{L-BD:KSepsm}) and
% the second bound in (\ref{L-BD:ChSh04bd})
(\ref{090615a})
yields
% for the second inequality,
%
%e3.7 ###
%
\begin{eqnarray}\label{L-BD:sbbound}\qquad
&&|\mathbb E [ h_{z,\varepsilon}(W)- \mathbb Eh_{z,\varepsilon}
(Z) ]|
\nonumber\\
&&\qquad\leq\frac{\mu}{\sigma^2}\Delta
%
%&\le&
%
%&\le&
+ \frac{\mu}{\sigma}\mathbb E
\biggl[\int_{(W^s-W)\wedge0}^{(W^s-W)\vee0}
(|W|+1)|t| \,dt
\nonumber\\
&&\qquad\quad\hspace*{57.2pt}{} +
\int_{-B/\sigma}^{B/\sigma}
\varepsilon^{-1}\int_{t \wedge0}^{t \vee0}
\mathbf{1}_{[z , z+\varepsilon]}(W+u) \,du\,dt \biggr]\\
&&\qquad\le\frac{\mu}{\sigma^2}\Delta+
\frac{\mu B^2}{2 \sigma^3}(\mathbb E|W|+1) +\frac{\mu}{\sigma
}\varepsilon
^{-1}\int_{-B/\sigma}^{B/\sigma} \int_{t \wedge0}^{t \vee
0}(0.4\varepsilon+2 D_Y)\,du\,dt\nonumber\\
&&\qquad\le
\frac{\mu}{\sigma^2}\Delta+ 1.4\frac{\mu}{\sigma^3}B^2 +\frac
{2\mu}{\sigma^3}B^2\varepsilon^{-1}D_Y,\nonumber
\end{eqnarray}
where in the second-to-last inequality above we have used the fact that
\[
P [\alpha\le W \le\beta] \le(\beta-\alpha)/\sqrt{2 \pi} + 2 D_Y,
\]
and in the last, the fact that $\mathbb E|W| \le1$.
By (\ref{090528a})
%Taking supremum over
%$\mathcaligr{H}$
%$z$
we see that $D_Y^\varepsilon$ is bounded by (\ref{L-BD:sbbound}),
and since
% by now noting from
(\ref{090528a}) and (\ref{090528b}) imply
%that
$
D_Y \le0.4\varepsilon+ D_Y^\varepsilon$,
%
% \mbox{substitution yields}
substitution yields
\[
D_Y \le\gamma(\varepsilon) :=\frac{a\varepsilon
+b}{1-c/\varepsilon},
\]
where
\[
a:=\frac{2}{5},\qquad
b:= \frac{\mu}{\sigma^2}\Delta+\frac{7}{5}\frac{\mu}{\sigma
^3}B^2\quad \mbox{and}\quad c:= \frac{2\mu B^2}{\sigma^3}.
\]
The optimum bound on $D_Y$ is at the positive root of
$\gamma'(\varepsilon)=0$, namely
$\varepsilon=c + r $ where
$r :=\sqrt{c^2+cb/a}$.

We wish to calculate $\gamma(c+r)$.
The denominator equals
\[
1-\frac{c}{c+r}=1-\frac{c(c-r)}{c^2-r^2}
%=1-\frac{c(c-r)}{-cb/a}
=1+\frac{a(c-r)}{b}=\frac{b+a(c-r)}{b},
\]
and therefore
% Now substituting this form of the denominator back into the function,
% we obtain
%
\begin{eqnarray*}
\gamma(c+r)
&=&
b \biggl( \frac{a(c+r)+b}{a(c-r)+b} \biggr)
%b ( \frac{c+r+b/a}{c-r+b/a} )
= b \biggl( \frac{c+b/a+r}{c+b/a-r} \biggr)
\\ &=&
b \biggl( \frac{(c+b/a+r)^2}{(c+b/a)^2-r^2} \biggr)
= b \biggl( \frac{(c+b/a+r)^2}{cb/a+(b/a)^2} \biggr)\\
&=& a \biggl( \frac{(c+b/a+r)^2}{c+(b/a)} \biggr)
%Now noting
%r=\sqrt{c}\sqrt{c+b/a}
%we obtain
% & = &
=
a \biggl( \frac{(c+b/a+\sqrt{c}\sqrt{c+b/a})^2}{c+(b/a)} \biggr)
\\
&= & a \bigl( \sqrt{c+b/a}+\sqrt{c} \bigr)^2
%
%jjjj
%Substituting the positive root
%now yields (\ref{090520a}) through simplifying
% ( \sqrt{c+b/a}+\sqrt{c} )^2\\
= \frac{2}{5} \Biggl(\sqrt{\frac{11}{2}\frac{\mu B^2}{\sigma
^3}+\frac{5}{2}\frac{\mu}{\sigma^2}\Delta} + \sqrt{\frac{2\mu
B^2}{\sigma^3}} \Biggr)^2,
\end{eqnarray*}
and this bound on $D_Y$
%now
yields (\ref{090520a}).
\end{pf}

Let $\operatorname{Bin}(n,p)$ denote the binomial distribution with
parameters $n \in\mathbb{N}$ and $p \in(0,1)$.
Our next two lemmas are concerned with binomial
and conditioned binomial distributions.
Lemma \ref{binlem} is used to prove Lemma
\ref{bincouplem}.
\begin{lemm}
\label{binlem}
Let $m \in\mathbb{N}$ and $p \in(0,1)$.
Suppose $N \sim\operatorname{Bin}(m,p)$,
and $\mathcaligr{L}(N') = \mathcaligr{L}(N | N >0)$,
$N''-1 \sim\operatorname{Bin}(m-1,p)$. Then
for all $k \in\mathbb{N}$,
%
%e3.8 ###
%
\begin{equation}\label{binlemeq}
P[ N \geq k ] \leq P [ N' \geq k ] \leq P[ N'' \geq k].
\end{equation}
\end{lemm}
\begin{pf}
The first inequality in (\ref{binlemeq}) is easy since
for $k \geq1$, by definition
$P[N' \geq k]
= P[N \geq k]/P[N \geq1]$.
It remains to prove the second inequality.
Suppose $\xi_1,\xi_2,\ldots$ are independent Bernoulli random
variables with parameter $p$. Let $M= \min\{ i\dvtx \xi_i =1\}$ and
$\tilde{N}'' := \sum_{i=M}^{M+m-1} \xi_i$.
Then $M$ and $\tilde{N}''$
are independent and $\mathcaligr{L}(\tilde{N}'') = \mathcaligr{L}(N'')$.

Define the random variables
\[
J:= \biggl\lceil\frac{M}{m} \biggr\rceil\quad\mbox{and}\quad
\tilde{N}' := \sum_{i=m(J-1)+1}^{mJ} \xi_{i}.
\]
In other words,
split the sequence of Bernoulli trials into disjoint intervals
of length~$m$, and let
$\tilde{N}'$ denote the number
of successful Bernoulli trials in the
first such interval that contains at least one successful trial.

Then $\tilde{N}'$ has the distribution of $N'$, and by construction
$ \tilde{N}' \leq\tilde{N}''$ almost surely.
Since $\tilde{N}''$ has the distribution of $N''$, this shows that
$N'$ is stochastically
dominated by $N''$, that is, the second inequality in (\ref{binlemeq}) holds.
\end{pf}

Our next lemma demonstrates the existence of a ``uniformly close
coupling'' of
random variables with a binomial distribution, and with
the same distribution conditioned to be nonzero (denoted,
respectively, $N$ and $M$ in the lemma).
This result will be used in Section \ref{couplsec} to
provide a uniformly close coupling of
$V$ [given by (\ref{Yndef})] and its size
biased version,
%via Lemma \ref{condsblem},
and likewise for
$S$ (in fact, for $n -S$).
%via Lemma \ref{sblem}.
%
\begin{lemm}
\label{bincouplem} Let $m \in\mathbb{N}$ and $p \in(0,1)$. Suppose $N
\sim\operatorname{Bin}(m,p)$, with
$N= \sum_{i=1}^m \xi_i$ where $\xi_i$ are independent
Bernoulli variables with parameter $p$.
Defining $\pi_k$
by
%
%e3.9 ###
%
\begin{equation}\label{pikdef}
\pi_k :=
\cases{
\dfrac{P[N > k |N >0]
- P[N > k] }{P[N=k](1 - (k/m))},
&\quad if $0 \leq k \leq m-1$,
\cr
0, &\quad if $k= m$,}
\end{equation}
we then have $0 \leq\pi_k \leq1$ for all $k \in\{0,\ldots,m\}$.

Suppose also that $\mathcaligr{B}$ is a further Bernoulli variable with
$P(\mathcaligr{B}=1|\xi_1,\ldots,\xi_m)=\pi_N$, and suppose $I$ is
an independent
discrete uniform random variable over $\{1,2,\ldots,m\}$.
Set $ M := N + (1-\xi_I)\mathcaligr{B}$,
that is, let $M$ be given by the same sum as $N$ except that
if $\mathcaligr{B}=1$ the $I$th term is set to 1.
Then
%
%e3.10 ###
%
\begin{equation}\label{1210a}
\mathcaligr{L}(M)=\mathcaligr{L}(N|N>0).
\end{equation}
\end{lemm}
\begin{pf}
Lemma \ref{binlem} shows $\pi_k \ge0$.
For the upper bound,
set $N'' = 1 + \sum_{i=2}^m \xi_i$.
Then $N'' - 1 \sim\operatorname{Bin}(m-1,p)$ and
$N'' $ is equal either to $N$ or to $N+1$, with
$P[N''=k+1 |N=k] = 1-(k/m)$
for $0 \leq k \leq m$.
Hence for all $k$,
by Lem\-ma~\ref{binlem},
\[
P[N > k] + P[N=k](1 -k/m) = P[N'' > k] \geq P[N > k |N > 0]
\]
so $\pi_k \leq1$.
%The final assertion, that
%$\mathcaligr{L}(M)=\mathcaligr{L}(N|N>0)$,
Also, assertion (\ref{1210a})
follows
by (\ref{pikdef}) and the fact that
\[
\{M > k \} = \{N > k\} \cup\{N=k, \mathcaligr{B}=1,
\xi_I=0\}.
\]
\upqed\end{pf}

Our next result refers to measurable real-valued functions
$\psi$ defined on all pairs $(x,\mathcaligr{X})$ such that
$\mathcaligr{X}$ is a finite subset of $C_n$ and $x \in\mathcaligr{X}$.
We say that such a functional $\psi$
is \textit{translation-invariant} if
$\psi(x,\mathcaligr{X}) = \psi(y+x,y+\mathcaligr{X})$
for all $x, \mathcaligr{X}$ and all $y \in C_n$ (here
addition is in the torus $C_n$, and $y + \mathcaligr{X}:= \{y+w \dvtx w
\in
\mathcaligr{X}\}$). For $ r >0$, we say that $\psi$ has \textit{radius
$r$} if
$\psi(x,\mathcaligr{X})$ is unaffected by the addition of points to,
or removal of
points from, the point set $\mathcaligr{X}$ at a distance more than
$r$ from $x$,
that is, if for all $(x,\mathcaligr{X})$ we have $\psi(x,\mathcaligr
{X}) = \psi(x,\mathcaligr{X}\cap B_r(x))$.
The notion of radius is the same as that of \textit{range of interaction}
used in \cite{PR}; see
also the notion of \textit{radius of stabilization}, in
\cite{PYStein,PR} and elsewhere. We also define
\begin{eqnarray*}
\| \psi\| &: =& \mathop{\operatorname{ess}\operatorname{sup}}_{x,\mathcaligr
{X}} \{ |\psi(x,\mathcaligr{X}) |\} ; \\
\operatorname{rng}( \psi) &: =&
\mathop{\operatorname{ess}\operatorname{sup}}_{x,\mathcaligr{X}} \{ \psi
(x,\mathcaligr{X}) \}
- \mathop{\operatorname{ess}\operatorname{inf}}_{x,\mathcaligr{X}} \{ \psi
(x,\mathcaligr{X}) \}.
\end{eqnarray*}
Recall that $\mathcaligr{U}_n := \{U_1,\ldots,U_n\}$
denotes a collection of $n$ independent
uniformly distributed points in $C_n$, and $\pi_d$ is
the volume of the unit $d$-ball.
\begin{lemm}
\label{momlem}
Let $n \in\mathbb{N}$ and $k \in\mathbb{N}$ with $2 \leq k \leq n$.
Suppose that for $i=1,\ldots, k$, $\psi_{i}$
is a measurable real-valued function
defined on all pairs $(x,\mathcaligr{X})$ with
$\mathcaligr{X}$ a finite set in $C_n$ and $x \in\mathcaligr{X}$.
Suppose for each $i$ that $\psi_{i}$ is translation-invariant
and has radius $r_i$ for some $r_i \in(0,\infty)$, and that
$\| \psi_{i} \| < \infty$,
and $\mathbb E[\psi_{1}(U_1,\mathcaligr{U}_n) ]=0$.
With $\phi_i := \pi_dr_i^d$,
suppose also that
$\phi_2 + \cdots+ \phi_k <n$.
Then
\begin{eqnarray*}
\Biggl| \mathbb E \Biggl[ \prod_{i=1}^k \psi_{i}(U_i,\mathcaligr
{U}_n) \Biggr] \Biggr|
&\leq&
\Biggl( n^{-1} \prod_{i=2}^k \|\psi_{i} \| \Biggr)
\operatorname{rng}(\psi_1)
\\
&&{}\times\Biggl( \pi_d
\Biggl( \sum_{i=2}^k (r_1 + r_i)^d \Biggr)
\\
&&\hspace*{17.8pt}{}
+ \phi_1 \Biggl( k-1 + \Biggl(\sum_{i=2}^k \phi_i \Biggr)
\Biggl(\frac{ 2n - \sum_{i=2}^k \phi_i }{n - \sum_{i=2}^k \phi_i }\Biggr)
\Biggr)\Biggr).
\end{eqnarray*}
\end{lemm}
\begin{pf}
Given $\mathbf{x}= (x_1,\ldots,x_k) \in C_n^k$, define the set of
points
\[
\mathcaligr{U}_n^\mathbf{x}:= \{x_1,\ldots,x_k, U_{k+1},\ldots, U_n
\}.
\]
Let $F_n$ be the set of $\mathbf{x}= (x_1,\ldots,x_k) \in C_n^k$
such that
$D(x_1,x_i) > r_1 + r_i $
for $i \in\{2,\ldots,k\}$, and let $F_n^c:= C_n^k \setminus F_n$.
Then by the law of total probability,
\begin{eqnarray*}
\mathbb E \Biggl[
\prod_{i=1}^k \psi_{i}(U_i,\mathcaligr{U}_n)
\Biggr]
&=& n^{-k} \int_{F_n} \mathbb E\prod_{i=1}^k \psi_{i}(x_i,\mathcaligr
{U}_n^\mathbf{x}) \,d\mathbf{x}
\\
&&{}
+ n^{-k} \int_{F_n^c} \mathbb E\prod_{i=1}^k \psi
_{i}(x_i,\mathcaligr{U}_n^\mathbf{x}) \,d\mathbf{x}.
\end{eqnarray*}
Since
$\mathbb E[\psi_1(U_1,\mathcaligr{U}_n)]=0$ it follows that
$\| \psi_1\| \leq\operatorname{rng}(\psi_1)$, so that
%
%e3.11 ###
%
\begin{eqnarray}\label{0422c}\qquad\quad
\Biggl|n^{-k} \int_{F_n^c} \mathbb E\prod_{i=1}^k \psi_{i}(x_i,\mathcaligr
{U}_n^\mathbf{x}) \,d\mathbf{x}\Biggr|
& \leq&
\Biggl( \prod_{i=1}^k \|\psi_{i} \| \Biggr)
P [(U_1,\ldots,U_k) \in F_n^c]
\nonumber\\[-8pt]\\[-8pt]
& \leq&
\operatorname{rng}(\psi_1)
\Biggl( \prod_{i=2}^k \|\psi_{i} \| \Biggr) \sum_{ i =2}^k
\pi_d(r_1+r_i)^d/n.\nonumber
\end{eqnarray}
Fix $\mathbf{x}= (x_1,\ldots,x_k) \in F_n$.
For $m \in\mathbb{Z}_+$,
let $h_1(m):= \mathbb E\psi_{1}(x_1,\{x_1\} \cup\mathcaligr{Y}_m)$, where
$\mathcaligr{Y}_m$ denotes
a collection of $m$ uniformly distributed points in $B_{r_1}(x_1)$.
Let
$h_2(m):= \mathbb E\prod_{i=2}^k \psi_{i}(x_i,\{x_2,\ldots,x_k\}
\cup\mathcaligr{Y}'_m)$,
where
$\mathcaligr{Y}'_m$ denotes
a collection of $m$ uniformly distributed points in
$\bigcup_{i=2}^k
B_{r_i}(x_i)$.

If $N_1$ and $N_2$ denote the number of points of
$\{U_{k+1},\ldots,U_n\}$ in
$B_{r_1}(x_1)$ and in $\bigcup_{i=2}^k B_{r_i}(x_i)$,
respectively,
then the values of
$\psi_{1}(x_1,\mathcaligr{U}^\mathbf{x}_n) $ and of
$ \prod_{i=2}^k \psi_i(x_i, \mathcaligr{U}^\mathbf{x}_n)$
are conditionally
independent, given $(N_1,N_2)$, because the regions $B_{r_1}(x_1)$ and
$\bigcup_{i=2}^k B_{r_i}(x_i)$ are disjoint since we assume $\mathbf
{x}\in F_n$.
Hence, we assert that
%
%e3.12 ###
%
\begin{equation}\label{0422d}
\mathbb E \Biggl[ \prod_{i=1}^k \psi_{i}(x_i,\mathcaligr{U}^\mathbf{x}_n)\Biggr]
= \mathbb E [ h_{1}(N_1) h_2(N_2) ],
\end{equation}
where $(N_1,N_2,N_3)$ have the multinomial distribution
%
%e3.13 ###
%
\begin{equation}\label{mul}
( N_1,N_2,N_3 ) \sim\operatorname{Mult}
\biggl( n-k;\frac{a_1}{n},
\frac{a_2}{n},
\frac{a_3}{n}\biggr)
\end{equation}
with $a_1$
denoting the volume of a ball of radius $r_1$ in $C_n$
[so that $a_1 \leq\phi_1 = \pi_dr_1^d$ with equality if
$r_1 \leq(1/2) n^{1/d}$]
while $a_2$ is the volume of
$\bigcup_{i=2}^k B_{r_i}(x_i)$
in $C_n$ and $a_3 := n-a_1 -a_2$.
To verify (\ref{0422d}),
use the law of total probability to decompose
the left-hand side as a sum over possible values
of $(N_1,N_2)$.

Also, if $\tilde{N}_1 \sim\operatorname{Bin}(n-1,\frac{a_1}{n})$
then $\mathbb E[ h_{1}(\tilde{N}_1)] = 0$,
because of the assumption that
$\mathbb E[\psi_{1}(U_1,\mathcaligr{U}_n) ]=0$,
along with translation invariance; the
value of $\mathbb E[ h_{1}(\tilde{N}_1)] $ does not depend on $x_1$.

We give a coupling of $N_1$ to another random variable
$N'_1$ with the same distribution as $\tilde{N}_1$
that is independent of $N_2$,
for which we can give a useful bound on
$P[N_1 \neq N'_1]$.

Consider throwing a series of colored balls
so each ball can land in one of three urns, where the probability
of landing in urn $i$ is $a_i/n$ for $1 \leq i \leq3$.
First, throw $n-k$ white balls and let $N_1^*, N_2,N^*_3$ denote
the number of
white balls in urn $i$ for $i = 1,2,3$, respectively, that is, let
$(N^*_1,N_2,N^*_3)$ have
the $\operatorname{Mult}(n-k;\frac{a_1}{n},\frac{a_2}{n},\frac{a_3}{n})$
distribution. Now pick out
the $n-k- N_2 $ balls in urns 1 and 3, paint them red,
and throw them again; that is,
given the values of $N^*_1,N_2,N^*_3$ let
$N_1^r,N_2^r,N_3^r$ count the number of red balls in urns
$1,2,3$, respectively,
and so be nonnegative integer valued
variables such that
\[
\mathcaligr{L} ( ( N_1^r,N_2^r,N_3^r ) |N_1^*,N_2)
= \operatorname{Mult}
\biggl( n- k - N_2 ;\frac{a_1}{n}, \frac{a_2}{n} , \frac{a_3}{n}\biggr).
\]

Now take the $N_2^r$ red balls in urn $2$,
paint them blue, and throw them again but condition them to land
in urns 1 and 3 (or equivalently, throw each blue ball
again and again until
it avoids urn 2), so that
\[
\mathcaligr{L} ( ( N_1^b,N_3^b
)|N^*_1,N_2,N_1^r,N_2^r )
=\operatorname{Mult} \biggl(N_2^r ; \frac{a_1}{a_1+a_3},\frac{
a_3}{a_1+a_3} \biggr).
\]
Finally, throw $k-1+N_2$ green balls, making the total number
of green, red and blue balls
$n-1$, and record how many land in urn 1, so
\[
\mathcaligr{L}
(N_1^g
| N^*_1,N_2,N_1^r,N_2^r,N_1^b ) =
\operatorname{Bin}
\biggl( k-1 + N_2; \frac{a_1}{n} \biggr).
\]
Now set
\[
N_1 = N_1^r + N_1^b ,\qquad
N_3 = N_3^r + N_3^b\quad \mbox{and}\quad N'_1 = N_1^r + N_1^g.
\]
Then $(N_1,N_2,N_3) $ have the multinomial distribution given
by (\ref{mul}).
Also, $N'_1 \sim\operatorname{Bin}(n-1,\frac{a_1}{n})$
and $N'_1$ is independent of $N_2$.

Since $N'_1 = N_1 - N_1^b +N_1^g$,
we have that
\begin{eqnarray*}
P[ N_1 \neq N'_1] & \leq&
\mathbb E[N_1^g]
+
\mathbb E[N_1^b]
\leq
\frac{a_1}{n}
(k-1 + \mathbb EN_2 )
+ \biggl(
\frac{a_1}{a_1 +a_3} \biggr) \mathbb E[N_2^r]
\\
& \leq&
\frac{a_1}{n}
( k-1 + a_2 )
+
\biggl( \frac{a_1}{n - a_2} \biggr)
a_2
\end{eqnarray*}
so that
\[
\bigl|
\mathbb E\bigl[ h_{2}(N_2) \bigl( h_{1}(N_1) -
h_1(N'_1) \bigr) \bigr]
\bigr|
\leq\frac{a_1}{n} \biggl(
k-1 + a_2
+
\biggl(
\frac{
n a_2
}{n - a_2}
\biggr)
\biggr) \operatorname{rng}(\psi_1) \prod_{i=2}^k \|\psi_{i}\|
\]
and since $N'_1$ is independent of $N_2$ with
$N'_1 \sim\operatorname{Bin}(n-1,\frac{a_1}{n})$,
\[
\mathbb E[ h_{1}(N'_1) h_{2}(N_2) ] = 0,
\]
so by (\ref{0422d})
and the fact that $a_1 \leq\phi_1$ and $a_2 \leq\sum_{i=2}^k \phi
_i $
and the assumption that $ \sum_{i=2}^k \phi_i < n$,
\begin{eqnarray*}
&&\Biggl|
\mathbb E \Biggl[ \prod_{i=1}^k \psi_{i}(x_i,\mathcaligr{U}_n^{\mathbf{x}})
\Biggr] \Biggr|
\\
&&\qquad\leq\frac{a_1}{n} \biggl(
k-1+ a_2
\biggl( \frac{ 2n - a_2 }{ n - a_2}
\biggr) \biggr) \operatorname{rng}(\psi_1) \prod_{i=2}^k \|\psi
_{i}\|
\\
&&\qquad\leq\frac{\phi_1}{n} \Biggl(
k-1+
\Biggl(
\sum_{i=2}^k \phi_i
\Biggr)
\biggl( \frac{ 2n -
\sum_{i=2}^k \phi_i
}{n -
\sum_{i=2}^k \phi_i
}
\biggr) \Biggr) \operatorname{rng}(\psi_1) \prod_{i=2}^k \|\psi
_{i}\|.
\end{eqnarray*}
The preceding bound holds uniformly over all possible values
of ${\mathbf x}= (x_1,\ldots,x_k) \in F_n$. Combined with (\ref{0422c}),
this shows that the asserted bound
holds.
\end{pf}
\begin{lemm}
\label{varbdlem}
Suppose $ \psi_1$ is as defined in Lemma \ref{momlem}.
Then with notation from that result, if $\phi_1 < n$ then
\begin{eqnarray*}
&&\operatorname{Var} \Biggl[ \frac{1}{n} \sum_{i=1}^n \psi
_1(U_i,\mathcaligr{U}_n) \Biggr]
\\
&&\qquad\leq
\frac{ \| \psi_1\|^2 }{n}
( 1 +
2^d
\phi_1
)
+
\frac{\|\psi_1\|
}{n}
% \frac{\|\psi_1\| \phi_1}{n}
\biggl(\phi_1+\phi^2_1\biggl( \frac{2 n -\phi_1 }{n -\phi_1} \biggr)
\biggr) \operatorname{rng}(\psi_1).
\end{eqnarray*}
\end{lemm}
\begin{pf}
By the case $k=2$ of Lemma \ref{momlem},
\begin{eqnarray*}
&&\operatorname{Cov}( \psi_1(U_1,\mathcaligr{U}_n),\psi
_1(U_2,\mathcaligr{U}_n))
\\
&&\qquad
= \mathbb E[ \psi_1(U_1,\mathcaligr{U}_n)
\psi_1(U_2,\mathcaligr{U}_n) ]
\\
&&\qquad
\leq\frac{\|\psi_1\|}{n}\biggl(2^d \phi_1\|\psi_1\|+\phi_1
\biggl(1 + \phi_1\biggl( \frac{2n- \phi_1}{n - \phi_1} \biggr)\biggr)
\operatorname{rng}(\psi_1)\biggr)
\end{eqnarray*}
and since
\begin{eqnarray*}
\operatorname{Var} \Biggl[ \frac{1}{n} \sum_{i=1}^n \psi
_1(U_i,\mathcaligr{U}_n) \Biggr]
&=&
n^{-1}
\operatorname{Var}[ \psi_1(U_1,\mathcaligr{U}_n) ]
\\
&&{}+ \frac{n-1}{n}
\operatorname{Cov}( \psi_1(U_1,\mathcaligr{U}_n),\psi
_1(U_2,\mathcaligr{U}_n)),
\end{eqnarray*}
the result follows.
\end{pf}

%s4 ###
\section{Size-biased coupling constructions}
\label{couplsec}

We now give a
simple
lemma which shows how to size-bias a random variable
%$Y$
that can be expressed as a conditional probability
of an event arising from some further randomization.
% beyond
%that of $Y$.
%is concerned with the construction
%of variables with size-biased distributions.
%
\begin{lemm}
\label{condsblem}
Suppose $Y$ is a random variable given by
$Y= a P[A|\mathcaligr{F}]$, where $\mathcaligr{F}$ is some $\sigma$-algebra,
$a>0$ is a constant,
%$aZ$ is a nondegenerate Bernoulli random variable
and
$A$ is an event with $0<P[A]< 1$.
Then
%a random variable
$Y'$ has the $Y$ size biased distribution
if
%
%e4.1 ###
%
\begin{equation}\label{newsbeq}
\mathcaligr{L}( Y') = \mathcaligr{L}(Y|A).
\end{equation}
\end{lemm}
\begin{pf}
With $\mathcaligr{L}(Y')$ defined by (\ref{newsbeq}),
we must show for all bounded and continuous $g\dvtx\mathbb{R}\to\mathbb{R}$,
% For such $g$,
that $\mathbb E[ g(Y') ] = \mathbb E[Y g(Y)]/ \mathbb E[Y]$
[see (\ref{formalsb})]. But
\begin{eqnarray*}
\mathbb E[g(Y')] & = &
\mathbb E[g(Y)|A] = \mathbb E[ g(Y)\mathbf{1}_A]/P[A]
\\
& = & \mathbb E[ g(Y) P[A|\mathcaligr{F}] ] / P[A],
\end{eqnarray*}
where the last equality follows because $g(Y)$ is $\mathcaligr{F}$-measurable.
The last expression equals $\mathbb E[Y g(Y)]/\mathbb E[Y]$, as
required.
\end{pf}

Let $V$ and $S$ be given by (\ref{Vndef}), (\ref{Yndef}), respectively.
% \eq{Yndef}. , $S$ be given by \eq{Vndef}, \eq{Yndef}.
Set $W=n-S$ (the number of nonsingletons).
%Then $V$
We assert that either $V$ or $W$
can be expressed as $n$ times the
conditional probability of some event $A$, given
the locations of the points of $\mathcaligr{U}_n$, so that
Lemma~\ref{condsblem} is applicable.
For $V$, take $A= A_V$ to be the event that
an additional uniformly distributed random point $U_0$ in $C_n$ lies
in the covered region $\bigcup_{i=1}^n B_{i,\rho}$.
For $W$, take $A = A_W$
to be the event that
an element of $\mathcaligr{U}_n$, selected uniformly at random,
is
% a non-singleton.
nonisolated.

%In the case of $V$,
%both cases (and this is why it helps to work with
%$W$ rather than $S$),
Event $A_V$ can be written as
the event that $N_V > 0$, where
% for a certain binomial
%$\Bin(n,n^{-1}\pi_d \rho^d)$
%variable which we denote by $N$.
% in either case, and
%where $m= n$ in the case of $V$ and $ m = n-1$ in the
%case of $W$.
%In the case of $V$,
$N_V$ denotes the
number of points of $\mathcaligr{U}_n$ in
$B_\rho(U_0)$, and
$N_V \sim\operatorname{Bin}(n,\phi/n)$ (recall $\phi:= \pi_d\rho^d$
and $C_n$ has volume $n$).
%where we set $m =n$.
%In the case of $W$, $N$ is the number of points
%of $\U_n \setminus\{U_0\}$ in $B_\rho(U_0)$ if we label the
%distinguished
%point of $\U_n$ as $U_0$.
%In either case,
A point set
(denoted $\mathcaligr{U}_V$)
with the
conditional distribution of $\mathcaligr{U}_n$ given $N_V$
%% (and also denoted $\U_n$) can be obtained as follows:
%(denoted $\U_F$ in the $V$ case and denoted $\U_W$ in the $W$ case)
can be obtained as follows:

\begin{enumerate}[III.]
\item[I.]
Sample a uniform random point in $C_n$, denoted $U_0$.
%(in the second case only, this will be part of $\U$)
%
\item[II.]
Set $m=n$.
Sample $N=N_V$ independent uniform random points in $B_\rho(U_0)$,
%the $\rho$-neighborhood of $U_0$, and
%Sample
and $m-N$ independent uniform random points in $C_n \setminus B_\rho(U_0)$.
\item[III.]
%In the case of $V$, let
Let $\mathcaligr{U}_V$ be the union of the two samples of uniform points.
%In the case of $W$, let
% $\U_W$ be the union of the two samples of uniform ramdom points,
%together with an added point at $U_0$.
\end{enumerate}
Therefore, coupled realizations of $\mathcaligr{U}_V$ and $\mathcaligr{U}'_V$
(having, respectively, the distribution of $\mathcaligr{U}_n$ and
the conditional distribution of $\mathcaligr{U}_n$ given $N_V>0$),
and hence coupled realizations of $V$ and $V'$,
can be obtained
%in the case of $V$
%by the following steps.
as follows.
% may be obtained
%in both cases as follows.
% Since the couplings are rather similar in the two cases,
%we describe them together, indicating those places where they
%diverge.

\begin{enumerate}
\item Set $m=n$.
%Sample $N$ from the $\Bin(m, \phi/n)$ distribution.
%with $m=n$ in the first case, or $m=n-1$ in the second case
%(cf. Lemma \ref{bincouplem})
%
\item
Sample $U_0$ uniformly at random over $C_n$.
\item
Sample $m$ random $d$-vectors
independently and uniformly over
%$B_\rho(U_0)$,
$C_n$,
and denote this
point set by $\mathcaligr{U}_{m,1}$.
\item Let $N$ denote the number of points of $\mathcaligr{U}_{m,1}$ in
$B_\rho(U_0)$.
\item
Sample a Bernoulli random variable
$\mathcaligr{B}$ with $P[\mathcaligr{B}=1] = \pi_N$,
where $(\pi_k, k \geq0)$ is given by (\ref{pikdef}).
%Sample $m-N$ random $d$-vectors
% independently and uniformly over $C_n \setminus B_\rho(U_0)$,
%and denote this point set by $\U_{m,2}$.
%Let $\U_V := \U_{n,1} \cup\U_{n,2}$.
%
\item Sample a random $d$-vector $U$ which is uniform over
$B_\rho(U_0)$.
\item
If $\mathcaligr{B}=1$, then select one of the points of
% $\U_{n,1} \cup\U_{n,2}$
$\mathcaligr{U}_{m,1} $
uniformly at random, and move it to $U$. Denote the resulting
modification of
%$\U_{n,1} \cup\U_{n,2}$
%by $\U_{n,3}$.
$\mathcaligr{U}_{m,1} $
by $\mathcaligr{U}_{m,2}$.
%If $\BB=0$ then set $\U_{n,3} = \U_{n,1} \cup\U_{n,2}$.
If $\mathcaligr{B}=0$ then set $\mathcaligr{U}_{m,2} : = \mathcaligr
{U}_{m,1} $.
\item
Set $\mathcaligr{U}_V := \mathcaligr{U}_{m,1} $ and
%Set
$\mathcaligr{U}'_V := \mathcaligr{U}_{m,2}$.
Set $V := g_V(\mathcaligr{U}_V)$ and
$V':=g_V(\mathcaligr{U}'_V)$, where $g_V(\mathcaligr{U}):=
\operatorname{Vol}( \bigcup_{x \in\mathcaligr{U}} B_{\rho}(x))$.
\end{enumerate}
By Lemma \ref{bincouplem},
the number of points of
$\mathcaligr{U}_{m,2}$ in $\mathcaligr{B}_\rho(U_0)$ has the distribution
$\mathcaligr{L}(N_V|N_V>0)$, and hence $\mathcaligr{L}(\mathcaligr
{U}'_V) = \mathcaligr{L}(\mathcaligr{U}_V|N_V >0)$.
So by Lemma \ref{condsblem}, $V'$ has the $V$ size biased distribution.

In the case of $W$, $A_W$ is the event that $N_W>0$,
where $N_W$ denotes the number of points of $\mathcaligr{U}_n
\setminus\{U_0\}$
in $B_\rho(U_0)$, and now $U_0$ denotes a
point of $\mathcaligr{U}_n$ selected uniformly at random.
So $N_W \sim\operatorname{Bin}(n-1,\phi/n)$. We can obtain
a point set (denoted $\mathcaligr{U}_W$)
with the conditional distribution of $\mathcaligr{U}_n$ given $N_W$
by the same steps as for $\mathcaligr{U}_V$ except that now in Step II
we put
$m=n-1$ and $N=N_W$, and in Step III, $\mathcaligr{U}_W$ is the union
of the two
samples of uniform random points with an added point at $U_0$.
Hence,
we can obtain coupled realizations
of
$W$ and $W'$
% $\U_W$ and $\U'_W$
by the same
sequence of steps as described above for $(V,V')$, except
that the following steps are modified:
\begin{itemize}
\item In Step 1, set $m= n-1$ (this affects Steps 3 and 5.)
\item
In Step 8,
%set $\U_W := \U_{n,1} \cup\U_{n,2} \cup\{U_0\}$.
set $\mathcaligr{U}_W := \mathcaligr{U}_{m,1} \cup\{U_0\}$,
%In Step 9, if $\BB=1$ then select a point of $\U_{n,1}
%of the resulting modification of $\U_{n,1} \cup\U_{n,2}$.
%In Step 10,
%and set $\U'_W := \U_{n,3} \cup\{U_0\}$,
and $\mathcaligr{U}'_W := \mathcaligr{U}_{m,2} \cup\{U_0\}$.
%In Step 11,
%In the $W$ case,
Set $W := g_W(\mathcaligr{U}_W)$ and $W':= g_W(\mathcaligr{U}'_W)$
with $g_W(\mathcaligr{U}) := \sum_{x \in\mathcaligr{U}} \mathbf{1}
\{\mathcaligr{U}\cap B_{\rho}(x) \neq\{x\} \}$.
\end{itemize}

By a similar argument to the $V$ case,
$W'$ has the $W$ size biased distribution.
%We shall need also to consider the point processes $\U_W$ and
%$\U'_W$ in the proof of Theorem \ref{thm1}.

%s5 ###
\section{Proof of Theorem \protect\ref{thm2}}
\label{proof-2}

%We can describe the size-biased distribution
%of $V$ using Lemma \ref{condsblem}.
%To do this, consider the measurable space $(E,\mathcaligr{E})$
%where $E$ is the disjoint union of $C_n^n$ and $C_n$.
%For $x$ and $t $ in $ E$, set
% $h(x,t) $ to be $\mbox{Vol}(C_n)=n$ if $x = (x_1,\ldots,x_n) \in
%C_n^n$,
%and $t \in C_n$, and $t \in\cup_{i=1}^n B_\rho(x_i)$, and
%set $h(x,t)$ to be zero otherwise.
%Let $X$ be a $C_n^n$-valued random variable with the distribution
%of $\UU= (U_1,\ldots,U_n)$, and let $T$ be a uniformly distributed
%variable in $C_n$, independent of $X$.
%Then the volume of the union of the
%$\rho$-balls centered at the points of $X$ is the same as
%$\E[h(X,T)|X]$, i.e. $V \eqd\E[h(X,T)|X]$.
%
%Under the assumption made here that points are uniformly distributed
%over the torus $C_n$, the variable $\tilde{T}$ given by \eq{tilYdef}
%is also uniformly distributed over the torus, and is denoted
%$U_0$ below. Given
% $U_0$, the $C_n^n$-valued variable
%$\tilde{X}$ given by \eq{tilXdef} is uniformly distributed
% over those $(x_1,\ldots,x_n) \in C_n^n$ for
%which $U_0$ lies in the union of the $\rho$-balls centered
%at the points $x_i, 1\leq i \leq n$.
%

%By Lemma \ref{condsblem},
We couple $V'$ to
$V$ as described in Section \ref{couplsec}.
% explained at the end of the preceding section.
%with the size-biased distribution of $V$ can
%be obtained by first sampling a point $U_0$ uniformly over $C_n$,
%then sampling $(U'_1,\ldots,U'_n)$ under a uniform distribution
%over $C_n^n$ conditional on at least one point $U'_i$ lying distance
%at most $ \rho$ from $U_0$, and then setting
%V' : = \mathrm{Volume} (
% \cup_{i=1}^n B_{\rho} (U'_i).
% )
%Given $U_0$, the number of points of
%$\{U'_1,\ldots,U'_n\}$ in $B_\rho(U_0)$ has the distribution
%of a $\Bin(n, \vold\rho^d/\mbox{Vol}(C_n))$
% variable conditioned to be at least 1
%%(the condition $6^d \phi<n$ in the statement of
%Theorem \ref{thm2} precludes any possibility
%of reduced volume of $B_\rho(U_0)$ due to wrap-around
%effects). Thus
%we can obtain a random variable $V''$ that is coupled to $V$ and
%has the same distribution as $V'$, as follows.
%Sample random $d$-vectors $(U_0,U_1,\ldots,U_{n})$,
% independent and uniform over $C_n$.
%Set $\U:= \{U_1,\ldots, U_{n}\}$.
%Sample a value of $I$ from the discrete uniform distribution on
% $\{1,\ldots,n\}$, independent of $\U$.
%Let $N$ denote the number of points of $\U$ in $B_\rho(U_0)$.
%Sample a Bernoulli random variable
% $\BB$ with $P[\BB=1] = \pi_N$,
%%[LG: slight additions MP: OK]
%where $(\pi_k, k \geq0)$ is given by
%(cf. Lemma \ref{bincouplem}.)
%$B_\rho(U_0)$.
%$$
%U_i'' = \{ \begin{array}{ll}
%U & \mathrm{if} I=i \mathrm{and} \BB=1
%U_i & \mathrm{otherwise}
%$$
%Define $V = V_n$ by \eq{Vndef}
%and
%set
%V'' : = \mathrm{Volume} (
% \cup_{i=1}^n B_{\rho} (U''_i).
% )
Since $V'$ differs from $V$ through the moving of at most a
single point, clearly $|V'-V| \leq\pi_d\rho^d := \phi$.
% Hence we may take $B=\volrho$
%in the application of Lemma \ref{LGthm}.
%
Hence, by Lemma \ref{LGthm} with $B = \phi$, to prove Theorem \ref
{thm2} it suffices
to prove the following.
\begin{prop}
\label{convar2prop}
Under the assumptions of Theorem \ref{thm2},
$\operatorname{Var}(\mathbb E[V'-V|V]) \leq n^{-1}\eta_V(n,\rho)$,
where $\eta_V(n,\rho)$ is given by (\ref{etaVdef}).
\end{prop}
\begin{pf}
Let $\mathcaligr{G}$
be the $\sigma$-algebra generated by the point set
$\mathcaligr{U}_V$.
%Information about
%%$U_0$, $U$, $N$ or
%which points
%of $\mathcaligr{U}_V$ came from where is \textit{not} included in
%$\mathcaligr{G}$.
List the points of $\mathcaligr{U}_V$, in an order chosen uniformly at random,
as $U_1,\ldots,U_n$, and set $\mathbf{U}:= (U_1,\ldots,U_{n})$.
Then $V$ is $\mathcaligr{G}$-measurable.
The conditional variance
formula, with $X = \mathbb E[V'-V|\mathcaligr{G}]$, yields
\[
\operatorname{Var}(\mathbb E[V'-V|V]) = \operatorname{Var}( \mathbb
E[X|V]) \leq\operatorname{Var}(X) ,
\]
so it suffices to prove
%
%e5.1 ###
%
\begin{equation}\label{a0422a}
\operatorname{Var}(\mathbb E[V'-V|\mathcaligr{G}])
\leq n^{-1} \eta_V(n,\rho).
\end{equation}

For $x \in C_n$,
let $\xi_x$ denote the probability that $\mathcaligr{B}=1$, given
$\mathcaligr{U}_n$ and given that $U_0=x$,
that is, $\xi_x=\pi_{N_x}$, where $N_x$ denotes the number
of points of $\mathcaligr{U}_V$ in $B_\rho(x)$.
Let $R_{xj}$ denote the expectation (over $U$) of the increment in
the covered volume
if $U_j$ is moved to a uniform
randomly selected location $U$ in $B_\rho(x)$.
Note that
%[LG: slight changes MP: more slight changes]
for $x$ and $j$ fixed, $R_{xj}$ is determined by
%$\{U_i\}_{i=1}^n$.
$\mathbf{U}$.
Then, since both $U_0$ and $I$ are independent of $\mathcaligr{G}$,
\[
\mathbb E[V'-V|\mathcaligr{G}] = \frac{1}{n} \int_{C_n} \xi_x
\Biggl( \frac{1}{n} \sum_{j=1}^n R_{xj} \Biggr)\,
dx,
\]
where the first factor of $1/n$ comes from the probability density
of $U_0$, and the second arises as the probability that $I$ takes the
value $j$.

Let $H_x$ be the expectation (over $U$) of the increment in the
covered volume when a point is inserted into $\mathcaligr{U}_V$ at a uniform
random location $U\in B_\rho(x)$, and let $T_j$ be the increment in
the covered volume when point $U_j$ is removed from $\mathcaligr{U}$
%[LG: slight changes. MP: more slight changes]
(for fixed
$x$ and $j$, both $H_x$ and $T_j$ are determined by $\mathbf{U}$).
% $(U_1,\ldots,U_n)$).
If $U_j$ is far distant from $x$ then $R_{xj}= H_x + T_j$.
Set $ Q_{xj} := R_{xj} - H_x - T_j $, which is in fact
the expectation (over $U$) of the total volume
of the otherwise uncovered regions lying within
distance $\rho$ both of $U$ and of $U_j$
(such regions contribute to $T_j$ but not to $H_x$ or $R_{xj})$.
% [LG: moved verbal description
%which appeared after the display to here. MP. OK]
Then
%
%e5.2 ###
%
\begin{eqnarray}\label{a0422b}\quad
\mathbb E[V'-V|\mathcaligr{G}]
& = & \frac{1}{n^2} \int_{C_n}
\sum_{j =1}^n \xi_x(H_x+T_j +Q_{xj})\, dx
\nonumber\\[-8pt]\\[-8pt]
& = & \frac{1}{n} \int_{C_n} \xi_x \Biggl(H_x + \frac{1}{n} \sum
_{j =1
}^n Q_{xj} \Biggr)\,dx + \frac{1}{n^2} \int_{C_n} \sum_{j=1}^n \xi_x
T_j \,dx .\nonumber
%= \vold+ \frac{1}{n} \int_{C_n} \sum_{i=1}^n \xi_x
% (S'_x + \frac{1}{n} \sum_{j =1 }^n
%Q_{xj} )dx
\end{eqnarray}

Set $\phi:= \pi_d\rho^d$.
We have that
%[LG: slight changes. MP: OK]
$0 \leq H_x \leq\phi$, $ 0 \geq T_j \geq- \phi$ and
%in particular, $0 \leq Q_{xj} \leq3 \volrho$.
%Also,
if $D(x,U_j) > 3 \rho$ then
%$R_{xj} = H_x +T_j$ so
$Q_{xj}=0$.
Moreover,
if $D(x,U_j) > 3 \rho$ for all $j \in\{1,\ldots,n\}$, then
$ H_x = \phi$ and
if $D(x,U_j) > \rho$ for all $j \in\{1,\ldots,n\}$, then
$\xi_x =1$.
Finally, $Q_{xj} \geq0$ and
\[
0 \le H_x + \frac{1}{n} \sum_{j=1}^n Q_{xj} \leq
H_x + \sum_{j=1}^n Q_{xj} \leq
\phi.
\]
Hence setting
\[
\tau_x := \xi_x \Biggl(
H_x + \frac{1}{n} \sum_{j=1}^n Q_{xj} \Biggr) - \phi,
\]
we have that $ -\phi\leq\tau_x \leq0$,
% {\tt seems the bound could .
%e taken to be theta sub d; if xi =0 tau is minus theta, while if
%xi=1 and S and Q take their maximum values, tau is 2 theta - theta
%= theta}
and
$\tau_x$ is determined by the collection of points
of $\mathcaligr{U}_n$ within distance $3\rho$ of $x$, and $\tau_x
=0$ if there
are no
such points of $\mathcaligr{U}_n$. We can rewrite (\ref{a0422b}) as
\[
\mathbb E[V'-V|\mathcaligr{G}] = \phi+
\frac{1}{n} \int_{C_n} \tau_x \,dx + \frac{1}{n} \Biggl( \sum_{j=1}^n
T_j \Biggr)
+ \frac{1}{n^2} \int_{C_n} \sum_{j=1}^n ( \xi_x -1) T_j \,dx.
\]
Recalling that $B_r(x): = \{ y \in C_n \dvtx D(x,y) \leq r\}$, let
$\Gamma_{i,r}$ be the set of points
$y \in B_{r}(U_i)$
%{\tt B $y \in C_n \cap B_2(U_i)$ {\tt B
%sub 2 is defined as a subset of C sub n, so no need to take
%intersection}
such that $D(y,U_i) < D(y,U_j)$ for all $j \in
\{1,\ldots,n\}
\setminus\{i\}$ (i.e., the intersection of the $r$-ball around
$U_i$ and the Voronoi cell of $U_i$ relative to $\mathcaligr{U}_n$).
Set
\[
S'_i := \int_{\Gamma_{i,3\rho}} \tau_x\,dx,\qquad S''_i := \int
_{\Gamma_{i,\rho}}
(\xi_x -1)\,dx.
\]
Then
%
%e5.3 ###
%
\begin{eqnarray}\label{0423c}
\mathbb E[V'-V|\mathcaligr{G}] &=& \phi+ \Biggl( \frac{1}{n}
\sum_{i=1}^n S'_i \Biggr) + \Biggl( \frac{1}{n} \sum_{j=1}^n T_j
\Biggr) + \Biggl( \frac{1}{n^2} \sum_{i=1}^n S''_iT_i \Biggr)
\nonumber\\[-8pt]\\[-8pt]
&&{} + \biggl( \frac{1}{n^2} \sum_{(i,j)\dvtx i \neq j} S''_iT_j
\biggr),\nonumber
\end{eqnarray}
and if we put
$b=\mathbb ET_i$ (which
does not depend on $i$), we have
\[
\frac{1}{n^2} \sum_{(i,j)\dvtx i \neq j} S''_i T_j
=
\frac{1}{n^2} \biggl( \sum_{(i,j)\dvtx i \neq j} S''_i ( T_j -b) \biggr)
+ \frac{b(n-1)}{n^2} \Biggl( \sum_{i=1}^n S''_i \Biggr),
\]
so by (\ref{0423c}),
\begin{eqnarray*}
&&\mathbb E[V'-V|\mathcaligr{G}] \\
&&\qquad = \phi
+ \frac{1}{n^2} \biggl( \sum_{(i,j)\dvtx i \neq j} S''_i ( T_j -b) \biggr)
\\
&&\qquad\quad{} + \frac{1}{n}\sum_{i=1}^n \bigl( S'_i
+ T_i + \bigl( n^{-1} T_i + (1 -n^{-1})b \bigr) S''_i \bigr).
\end{eqnarray*}
Since $(x+y)^2 \leq2(x^2 +y^2)$ for any real $x,y$,
%
%e5.4 ###
%
\begin{eqnarray}\label{0525}\qquad
&&\operatorname{Var}(\mathbb E[V'-V|\mathcaligr{G}])
\nonumber\\
&&\qquad\leq
2\operatorname{Var} \Biggl( \frac{1}{n}\sum_{i=1}^n \bigl( S'_i
+ T_i + \bigl(n^{-1} T_i + (1 -n^{-1}) b\bigr) S''_i \bigr) \Biggr)
\\
&&\qquad\quad{} + 2 \operatorname{Var} \biggl(\frac{1}{n^2} \biggl(
\sum_{(i,j)\dvtx i \neq j} S''_i ( T_j -b) \biggr)\biggr).\nonumber
\end{eqnarray}

%We have the following table of radii and bounds on possible values (Table
Table \ref{table1} summarizes upper and lower bounds and the radius of
the relevant variables, where the \textit{radius} of a variable indexed by $i$ is
the smallest distance from $U_i$ one needs to look to establish its value
(as with the functionals considered in Lemma \ref{momlem}).

%
%t1
\begin{table}[b]
\caption{Radii and bounds for covered volume}
\label{table1}
\begin{tabular*}{\tablewidth}{@{\extracolsep{\fill}}lccccc@{}}
\hline
\textbf{Variable} & $\bolds{\tau_x}$ & $\bolds{T_i}$ & $\bolds{S''_i}$
& $\bolds{S'_i}$ & $\bolds{(n^{-1}T_i + (1-n^{-1})b)S''_i}$\\
\hline
Radius & $3 \rho$ & $2 \rho$ & $2\rho$ & $6 \rho$ & $2\rho$\\
Lower bound & $-\phi$ & $-\phi$ & $-\phi$ & $-3^d \phi^2$ & 0\\
Upper bound & 0 & 0 & 0 & 0 & $\phi^2$\\
\hline
\end{tabular*}
\legend{Note: The last two columns are deduced from the previous columns.}
\end{table}

Hence, the variable
\[
S'_i +
T_i + \bigl(n^{-1} T_i + (1 -n^{-1}) b\bigr) S''_i
\]
has radius $6\rho$ relative to $U_i$ and
lies between $-\phi- 3^d \phi^2$ and
$\phi^2$,
so that
its centered value is bounded
in absolute value by $(3^d+1)\phi^2 + \phi
$, and this also bounds its range of possible values.
So by Lemma \ref{varbdlem} and the assumption that $6^d \phi< n$,
%
%e5.5 ###
%
\begin{eqnarray}\label{0527c}\quad
&&\operatorname{Var} \Biggl( \frac{1}{n} \sum_{i=1}^n \bigl( S'_i + T_i + (n^{-1}
T_i + (1 -n^{-1}) b) S''_i \bigr) \Biggr)
\nonumber\\[-8pt]\\[-8pt]
&&\qquad\leq\frac{\phi^2 ( (3^{d} +1)\phi+ 1 )^2}{n} \biggl( 1 + (2^d +1) 6^d
\phi + \biggl( \frac{2 n - 6^d \phi }{ n - 6^d \phi } \biggr) 6^{2d} \phi^2 \biggr)
.\nonumber
\end{eqnarray}
Now consider the last term in the right-hand side of (\ref{0525}).
Set %$\bar{S}''_i := S''_i -a$ and
$\bar{T}_j:= T_j - b$. Then
%
%e5.6 ###
%
\begin{eqnarray}\label{0527a}\quad
&&\operatorname{Var} \biggl( \sum_{(i,j)\dvtx i \neq j} S''_i \bar{T}_j\biggr)
\nonumber\\
&&\qquad= n(n-1)(n-2)(n-3) \operatorname{Cov}( S''_1 \bar{T}_2,
S''_3\bar{T}_4)
\nonumber\\
&&\qquad\quad{} + n(n-1)(n-2) \bigl( \operatorname{Cov}(S''_1 \bar{T}_2,
S''_1 \bar{T}_3)
\\
&&\qquad\quad\hspace*{86.6pt}{}+ \operatorname{Cov}(S''_2 \bar{T}_1, S''_3 \bar{T}_1)
+ 2\operatorname{Cov}(S''_1 \bar{T}_2, S''_3 \bar{T}_1)\bigr)
\nonumber\\
&&\qquad\quad{} + n(n-1) \bigl(\operatorname{Var}( S''_1 \bar{T}_2 ) +
\operatorname{Cov}(S''_1 \bar{T}_2, S''_2 \bar{T}_1 )\bigr).\nonumber
\end{eqnarray}
It follows from the case $k=4$ of Lemma \ref{momlem} and
the assumption $6^d \phi< n$ (which implies $3(2^d \phi) < n$)
that
\begin{eqnarray*}
\operatorname{Cov}( S''_1 \bar{T_2}, S''_3 \bar{T_4} )
& = &
\mathbb E[S''_1 \bar{T_2}S''_3 \bar{T}_4 ] - (\mathbb E[S''_1 \bar{T}_2 ])^2
\leq \mathbb E[ S''_1 \bar{T_2}S''_3 \bar{T}_4 ]
\\
& \leq& \frac{\phi^4}{n} \biggl( 3\pi_d(4\rho)^d + 2^d \phi\biggl(3 + 3 (2^d) \phi
\biggl( \frac{2n- 3(2^d) \phi}{n- 3(2^d) \phi} \biggr) \biggr) \biggr)
\\
& = & \frac{3 \phi^4}{n} \biggl( (4^d + 2^d) \phi+ 4^d \phi^{2} \biggl(\frac{2n-
3(2^d) \phi}{n- 3(2^d) \phi}\biggr)\biggr).
\end{eqnarray*}
Since we can always bound
$\operatorname{Cov}(S''_i\bar{T}_j,S''_{i'} \bar{T}_{j'})$
above
by $\phi^4$, we have from (\ref{0527a}) that
%
%e5.7 ###
%
\begin{eqnarray}\label{0527b}
&& \operatorname{Var} \biggl( \frac{1}{n^2}
\biggl( \sum_{(i,j)\dvtx i \neq j} S''_i ( T_j -b) \biggr) \biggr)
\nonumber\\[-8pt]\\[-8pt]
&&\qquad \leq \frac{\phi^4}{n} \biggl( 3 (4^d + 2^d) \phi+ 3 (4^d) \phi^{2} \biggl(
\frac{2n- 3(2^d) \phi} {n- 3(2^d) \phi} \biggr) + 4 + \frac{2}{n}
\biggr).\nonumber
\end{eqnarray}

By (\ref{0525}), (\ref{0527c}) and (\ref{0527b}) we have that
\begin{eqnarray*}
&& (n/2) \operatorname{Var} ( \mathbb E[V'-V|\mathcaligr{G}] )
\\
&&\qquad \leq %3^{2d}
\phi^2 \bigl( (3^{d} +1)\phi+ 1\bigr)^2 \biggl( 1 + (2^d +1) 6^d \phi+ \biggl( \frac{2 n -
6^d \phi}{ n - 6^d \phi} \biggr) 6^{2d} \phi^2\biggr)
\\
&&\qquad\quad{} + \phi^4\biggl(3(4^d + 2^d) \phi+ 3(4^d) \phi^{2} \biggl( \frac{2n- 3(2^d) \phi}
{n- 3(2^d) \phi} \biggr) + 4 + \frac{2}{n} \biggr).
\end{eqnarray*}
This completes the proof of
Proposition \ref{convarprop}, and hence of Theorem \ref{thm2}.
\end{pf}

%s6 ###
\section{Proof of Theorem \protect\ref{thm1}}
\label{proof-1}

%Set $W=n-S$ (the number of non-singletons).
We couple $W'$ to
$W$ as described in Section \ref{couplsec}.
Thus $W= g_W(\mathcaligr{U}_W)$ and $W'= g_W(\mathcaligr{U}'_W)$, where
%$\UU'$
$\mathcaligr{U}'_W$ is obtained from
%$\UU$
$\mathcaligr{U}_W$ by moving at most a single
randomly selected point of $\mathcaligr{U}_W \setminus\{U_{0} \} $ to
a (uniform random) location in $B_\rho(U_0)$, if $\mathcaligr{B}=1$,
and leaving $\mathcaligr{U}_W$ unchanged if $\mathcaligr{B}=0$.

%Let $W$ be defined by \eq{Yndef} and set
%W'' : = \sum_{i=1}^n \mathbf{1} \{
%%\cap
%B_{\rho} (U''_i)
%Then
% $\LL(W'')= \LL(W')$.
%

The number of points that can be made isolated by removing a single
point from $\mathcaligr{U}_W$ is almost surely bounded by
$\kappa_d$. Moreover, the number of points that can be made
nonisolated by inserting a point (including the inserted
point itself) is almost surely bounded by $\kappa_d + 1$.
Hence $|W-W'| \leq\kappa_d + 1$, so we may
take $B=\kappa_d + 1$.

%Since $\Var(S)= \Var(W)$,
%it is easily seen that
%$(S - \E S)/\SD(S) = - (W - \E W)/\SD(W) $.
By the symmetry of the normal distribution,
%$d_K(\LL(X),Z)$ is the same as
%$d_K(\LL(-X),Z)$ for any random variable $X$, in particular for
%$X = (S - \E S)/\SD(S) $.
$D_{-S} = D_{S}$ and hence $D_W = D_S$.
Thus,
Theorem \ref{thm1}
follows from Lemma \ref{LGthm}
along with the following:
\begin{prop}
\label{convarprop}
Under the assumptions of Theorem \ref{thm1},
$\operatorname{Var}(\mathbb E[W'-W|W]) \leq
n^{-1}
\eta_S(n,\rho)$, where $\eta_S(n,\rho)$ is
given by (\ref{etaYdef}).
\end{prop}
\begin{pf}
Here we let $\mathcaligr{G}$ denote
the $\sigma$-algebra generated by the unlabelled point set
$\mathcaligr{U}:= \mathcaligr{U}_W$.
%, containing no information about which
%points of $\mathcaligr{U}$ come from where
Then $W'$ is $\mathcaligr{G}$-measurable,
and by the conditional variance formula
(as in the proof of Proposition \ref{convar2prop}),
it suffices to prove that
%
%e6.1 ###
%
\begin{equation}\label{0422a}
\operatorname{Var}(\mathbb E[W'-W|\mathcaligr{G}]) \leq n^{-1} \eta
_S(n,\rho).
\end{equation}

Label the points of $\mathcaligr{U}$, in an order chosen uniformly at random,
as $U_1,\ldots,U_n$, and set $\mathbf{U}:= (U_1,\ldots,U_n)$.
% $\xi_i$
%denote the conditional probability that $\BB=1$,
%given $\UU$ and given that $I=i$, i.e.
$\xi_i=\pi_{N_i}$, where
$N_i$ denotes the number of points of $\mathcaligr{U}\setminus\{U_i\}
$ in
$B_\rho(U_i)$. Let $R_{ij}$ denote the
expectation (over $U$) of the increment in the number of
nonisolated points when $U_j$ is moved to a uniform
randomly selected location $U$ in $B_\rho(U_i)$.
Then
\[
\mathbb E[W'-W|\mathcaligr{G}] = \frac{1}{n(n-1)} \sum_{(i,j)\dvtx i \neq
j} \xi_i R_{ij} ,
\]
where $\sum_{(i,j)\dvtx i \neq j}$ denotes summation over
pairs of distinct integers $i,j$ in $[1,n]$.

Now let $H_i$ be the expectation (over $U$)
of the increment in
the number of isolated points
when a point is inserted into $\mathcaligr{U}$ at a uniform random location
$U\in B_\rho(U_i)$, and let $T_j$ be the increment in
the number of isolated points
when point $U_j$ is removed from $\mathcaligr{U}$ (both $H_i$ and
$T_j$ are
determined by $\mathbf{U}$).
If $U_j$ is far distant from $U_i$ then $R_{ij}= - H_i - T_j$.
In fact, setting $ Q_{ij} := R_{ij} + H_i + T_j $,
we have that $Q_{ij}$ is
the expectation (over $U$) of the
number of otherwise isolated points of $\mathcaligr{U}$
within distance $\rho$ both of
$U$ and of $U_j$
(such points contribute to $T_j$ but not to $H_i$ or $R_{ij})$.
Then
%
%e6.2 ###
%
\begin{eqnarray}\label{0422b}
\mathbb E[W'-W|\mathcaligr{G}]
&=& \frac{1}{n(n-1)} \sum_{(i,j)\dvtx i \neq j} \xi_i(-H_i-T_j +Q_{ij})
\nonumber\\[-8pt]\\[-8pt]
&=& \frac{1}{n} \sum_{i=1}^n \xi_i \tau_i
- \frac{1}{n(n-1)} \sum_{(i,j)\dvtx i \neq j} \xi_i T_j,\nonumber
\end{eqnarray}
where we set
%
%e6.3 ###
%
\begin{equation}\label{taudef}
\tau_i := -H_i + \frac{1}{n-1} \sum_{j \dvtx j\neq i} Q_{ij}.
\end{equation}
Put $a := \mathbb E[\xi_i]$
(given $n$, this expectation does not
depend on $i$)
and put $b := (\kappa_d -1)/2$. Then
\begin{eqnarray*}
\frac{1}{n(n-1)} \sum_{(i,j)\dvtx i \neq j} \xi_i T_j &=&
\frac{1}{n(n-1)} \biggl( \sum_{(i,j)\dvtx i \neq j} (\xi_i -a) ( T_j - b)\biggr)
%T_j
\\
&&{} + \frac{a}{n} \Biggl( \sum_{j=1}^n
% ( T_j -b)
T_j \Biggr) + \frac{b}{n} \Biggl( \sum_{i=1}^n ( \xi_i -a)\Biggr).
%T_j
\end{eqnarray*}
Hence we can rewrite (\ref{0422b}) as
\[
\mathbb E[ W'-W|\mathcaligr{G}] = \frac{1}{n} \sum_{i=1}^n \bigl( \xi_i
\tau_i - aT_i - b (\xi_i-a) \bigr)- \frac{1}{n(n-1)} \sum_{(i,j)\dvtx i
\neq j} (\xi_i-a) (T_j-b) .
\]
Since $(x+y)^2 \leq2(x^2 + y^2)$ for any real $x,y$,
it follows that
%
%e6.4 ###
%
\begin{eqnarray}\label{0513c}\quad
\operatorname{Var} ( \mathbb E[W'-W|\mathcaligr{G}] ) &\leq& 2
\operatorname{Var} \Biggl( \frac{1}{n} \sum_{i=1}^n \bigl( \xi_i (\tau_i -b) +
a(b-T_i) \bigr) \Biggr)
\nonumber\\[-8pt]\\[-8pt]
&&{} + 2 \operatorname{Var} \biggl( \frac{1}{n(n-1)} \sum_{(i,j)\dvtx i \neq
j} (\xi_i-a) ( T_j -b)\biggr).\nonumber
\end{eqnarray}
We have that $-\kappa_d \leq H_i \leq0$, $ -1 \leq T_j \leq\kappa_d$,
and $Q_{ij} \geq0$. Also,
\[
0 \le-H_i + \sum_{j\dvtx j \neq i} Q_{ij} \leq\kappa_d,
\]
and if
$D(U_i,U_j) >3\rho$ then
$Q_{ij}=0$.
%
%t2
\begin{table}[b]
\caption{Radii and bounds for singletons}
\label{table2}
\tabcolsep=0pt
\noindent\begin{tabular*}{\tablewidth}{@{\extracolsep{\fill}}lccccccc@{}}
\hline
\textbf{Variable} & $\bolds{H_i}$ & $\bolds{T_i}$ & $\bolds{\xi_i}$
& $\bolds{\tau_i}$ & $\bolds{\xi_i(\tau_i-b)}$ &
$\bolds{a(b- T_i)}$ & $\bolds{\xi_i( \tau_i-b) + a(b- T_i)}$ \\
\hline
Radius & $3 \rho$ & $2\rho$ & $\rho$ & $3\rho$ & $3\rho$ &
$2 \rho$ & $3 \rho$\\
$\operatorname{ess}\operatorname{inf}$ & $-\kappa_d$ & $-1$ & 0 & 0
& %\frac{1 - \kappa_d}{2}
$(1 - \kappa_d)/2$ & $- ( \kappa_d +1 )/2$ & $- \kappa_d$\\
$\operatorname{ess}\operatorname{sup}$ & 0 & $\kappa_d$ & 1 & $\kappa_d$
& $( \kappa_d +1)/2$ & $( \kappa_d +1)/2$ & $\kappa_d + 1$\\
\hline
\end{tabular*}
\legend{Note: The last three columns are deduced from
the preceding columns and the definitions of $a,b$.}
\end{table}
Hence,
$ 0 \leq\tau_i \leq\kappa_d$, and
$\tau_i$ is determined by the collection of points
of $\mathcaligr{U}$ within distance $3 \rho$ of $U_i$.
%To summarize, if the \textit{radius} of a variable indexed by $i$
%is the smallest distance from $U_i$ one needs to look
%to establish its value (as with the functionals considered in
%Lemma \ref{momlem}),
%we have Table \ref{table2}.
Table \ref{table2} summarizes this discussion; recall from Table \ref{table1}
the notion of radius.

From the last column in this table, we see that after centering,
the terms in first sum in the right-hand side of (\ref{0513c})
have radius $3\rho$ and
absolute values bounded by
$ 1 + 2 \kappa_d $.
Moreover, even after centering each of these terms
has range (i.e., essential supremum minus essential infimum)
which is also bounded by $1 + 2 \kappa_d$ (this range is
unaffected by the centering).
Hence with $\phi:= \pi_d\rho^d$, Lemma \ref{varbdlem},
using the assumption
$3^d \phi< n$, yields
%
%e6.5 ###
%
\begin{eqnarray}\label{0519b}
&&\operatorname{Var} \Biggl( \frac{1}{n} \sum_{i=1}^n \bigl( \xi_i (\tau_i -b) + a(b
- T_i) \bigr) \Biggr)
\nonumber\\[-8pt]\\[-8pt]
&&\qquad\leq\frac{(1 + 2 \kappa_d)^2}{n} \biggl(1 + (2^d+1) 3^d \phi + \biggl( \frac{2n-
3^d \phi}{ n - 3^d \phi} \biggr) (3^d\phi)^{2}\biggr).\nonumber
\end{eqnarray}
Now consider the second sum in the right-hand side of (\ref{0513c}).
Set $\bar{\xi}_i := \xi_i -a$
and $\bar{T}_j:= T_j - b$.
Then
%
%e6.6 ###
%
\begin{eqnarray}\label{0519a}
&&\operatorname{Var} \biggl( \sum_{(i,j)\dvtx i \neq j} (\xi_i -a)( T_j -b)\biggr)
\nonumber\\
&&\qquad = n(n-1)(n-2)(n-3) \operatorname{Cov}( \bar{\xi}_1 \bar{T}_2,
\bar{\xi}_3\bar{T}_4)
\nonumber\\
&&\qquad\quad{} + n(n-1)(n-2) \bigl(
\operatorname{Cov}(\bar{\xi}_1 \bar{T}_2, \bar{\xi}_1 \bar{T}_3)\\
&&\qquad\quad\hspace*{86.06pt}{}
+ \operatorname{Cov}(\bar{\xi}_2 \bar{T}_1, \bar{\xi}_3 \bar{T}_1)
+ 2\operatorname{Cov}(\bar{\xi}_1 \bar{T}_2, \bar{\xi}_3 \bar{T}_1)\bigr)
\nonumber\\
&&\qquad\quad{} + n(n-1) \bigl(\operatorname{Var}( \bar{\xi}_1 \bar{T}_2 )
+ \operatorname{Cov}(\bar{\xi}_1 \bar{T}_2,
\bar{\xi}_2 \bar{T}_1)\bigr).\nonumber
\end{eqnarray}
Note that $\bar{\xi}_i$ has absolute value bounded by 1, and range
of possible values also bounded by 1, and mean zero.
Also, $\bar{T}_j$ has absolute value almost surely bounded by
$(\kappa_d + 1)/2$ (its mean
might not be zero). Hence,
the case $k=4$ of
Lemma \ref{momlem} [taking $r_1=r_2= \rho$
and $r_3=r_4 = 2 \rho$ so that $\phi_2+ \phi_3
+ \phi_4 = (2^{d+1} +1)\phi$] yields
\begin{eqnarray*}
&&\operatorname{Cov}( \bar{\xi}_1 \bar{T}_2, \bar{\xi}_3 \bar{T}_4 ) \\
&&\qquad=
\mathbb E[ \bar{\xi}_1 \bar{T}_2\bar{\xi}_3 \bar{T}_4 ] - (\mathbb
E[\bar{\xi_1} \bar{T}_2 ])^2 \leq \mathbb E[ \bar{\xi}_1
\bar{T}_2\bar{\xi}_3 \bar{T}_4 ]
\\
&&\qquad\leq\frac{ (\kappa_d + 1)^2 }{4 n} \biggl( \phi\bigl( 2(3^d) + 2^d\bigr) + 3 \phi+
(2^{d+1}+1) \phi^2 \biggl( \frac{2n - (2^{d+1} +1) \phi}{n - (2^{d+1} +1)
\phi} \biggr) \biggr),
\end{eqnarray*}
where we have also used the assumption
that $(2^{d+1} +1) \phi< n$.
Since we can always bound
$\operatorname{Cov}(\bar{\xi}_i\bar{T}_j,\bar{\xi}_{i'} \bar{T}_{j'})$
by $((\kappa_d +1)/2)^2$, we have from (\ref{0519a}) that
%
%e6.7 ###
%
\begin{eqnarray}\label{0519d}
&&\operatorname{Var} \biggl( \frac{1}{n(n-1)} \biggl(\sum_{(i,j)\dvtx i
\neq j} (\xi_i -a) ( T_j - b) \biggr) \biggr)
\nonumber\\
&&\qquad\leq\frac{
(\kappa_d +1)^2 }{4 n}
\biggl(\bigl( 2(3^d) + 2^d +3\bigr) \phi\nonumber\\[-8pt]\\[-8pt]
&&\hspace*{83pt}{} + (2^{d+1}+1) \phi^2 \biggl(
\frac{2n - (2^{d+1} +1) \phi }{n - (2^{d+1} +1) \phi} \biggr) \biggr)
\nonumber\\
&&\qquad\quad{} + \biggl(\frac{ \kappa_d +1}{2} \biggr)^2
\biggl(\frac{4}{n} + \frac{2}{n(n-1)} \biggr).\nonumber
\end{eqnarray}
By (\ref{0513c}), (\ref{0519b}) and (\ref{0519d}),
we have that
\begin{eqnarray}
&& n \operatorname{Var} ( \mathbb E[W'-W|\mathcaligr{G}] )
\nonumber\\
&&\qquad\leq 2 (1+ 2\kappa_d)^2 \biggl( 1 + (2^d+1) 3^d \phi+ \biggl(
\frac{2n- 3^d \phi}{ n - 3^d \phi} \biggr) 9^d\phi^2 \biggr)
\nonumber\\
&&\qquad\quad{} + \frac{(\kappa_d +1)^2}{2} \biggl( \bigl( 2(3^d) + 2^d
+3\bigr) \phi + (2^{d+1}+1) \biggl( \frac{2n - (2^{d+1} +1) \phi }{n -
(2^{d+1} +1) \phi} \biggr) \phi^{2} \biggr)
\nonumber\\
&&\qquad\quad{} + \frac{ (\kappa_d +1)^2}{2} \biggl( 4 + \frac{2}{n-1}
\biggr). \nonumber
\end{eqnarray}
This completes the proof of
Proposition \ref{convarprop}, and hence of Theorem \ref{thm1}.
\end{pf}

%s7 ###
\section{Proof of Theorem \protect\ref{thmlimsup} and numerics}
\label{secvar}

Again set $\phi:= \pi_d\rho^d$.
It is easy to see that provided $2 \rho< n^{1/d}$,
%
%e7.1 ###
%
\begin{equation}\label{0522a}
\mathbb E[V] = n \bigl(1-(1 - \phi/n)^{n}\bigr);\qquad
\mathbb E[S] = n (1 - \phi/n)^{n-1} ,
\end{equation}
and (\ref{meanlim}) follows from this.

Write $|\cdot|$ for the Euclidean norm and
recall that $\omega_d(|x|)$ denotes the volume of the union of unit
balls centered at the origin $\mathbf{0}$ and at $x$.
If $I_x$ denotes the indicator of the event that
$x$ is not contained in any of the balls $B_{\rho,i}$,
then
provided $4\rho< n^{1/d}$
we have the exact formula
%
%e7.2 ###
%
\begin{eqnarray}\label{varV}
\operatorname{Var}(V) &=& \operatorname{Var}(n-V)\nonumber\\
&=& \operatorname
{Var}\int_{C_n} I_x \,dx
\nonumber\\
&=& \int_{C_n} \int_{C_n} \mathbb E[I_x I_y] \,dx \,dy
- \bigl(n ( 1 - \phi/n)^n \bigr)^2
\\
&=& n\int_{B_{2\rho}(\mathbf{0})}
\biggl( 1 - \frac{\rho^d \omega_d(|y|/\rho) }{n} \biggr)^n \, dy
\nonumber\\
&&{}
+ n(n- 2^d\phi) \biggl( 1 - \frac{ 2\phi}{n} \biggr)^n
- n^2 ( 1 - \phi/n)^{2n}.\nonumber
\end{eqnarray}
\begin{pf*}{Proof of (\protect\ref{varVlim})}
For asymptotics as $n \to\infty$ with $\rho$ fixed,
use the MacLaurin expansion of $\log(1-x)$
to obtain
\begin{eqnarray*}
\biggl(1 - \frac{2\phi}{n} \biggr)^n &=& e^{-2 \phi} \exp\biggl( -
\frac{2 \phi^2}{n} + O( n^{-2} ) \biggr);
\\
\biggl(1 - \frac{\phi}{n} \biggr)^{2n} &=& e^{-2 \phi} \exp\biggl( -
\frac{\phi^2}{n} + O(n^{-2}) \biggr)
\end{eqnarray*}
so that
\begin{eqnarray*}
n^{-1} \operatorname{Var}( V ) & = & \int_{B_{2\rho}(\mathbf{0})}
\biggl( 1 - \frac{\rho^d \omega_d(|y|/\rho) }{n} \biggr)^n \, dy
\\
&&{} + ne^{-2\phi} \biggl( \biggl(1 - \frac{ 2^d \phi}{n} \biggr)
\exp\biggl( - \frac{2\phi^2}{n} \biggr)
- \exp\biggl( - \frac{\phi^2}{n} \biggr) + O(n^{-2}) \biggr)
\\
&\to& \biggl( \int_{B_{2\rho}(\mathbf{0})}
% ( 1 - \frac{\rho^d \omega_d(|y|/\rho) }{n} )^n
\exp\bigl( - \rho^d \omega_d(|y|/\rho) \bigr) \,dy \biggr) - e^{-2 \phi}
( 2^d \phi+ \phi^2)
\end{eqnarray*}
and this limit is equal to $g_V(\rho)$ as defined by
(\ref{gVdef}), so the first part of (\ref{varVlim}) is proven.

It remains to show that $g_V(\rho)>0$.
This can be done either by
using
the last part of Theorem 2.1 of
\cite{PY1}, or directly. We leave it
to the reader to check that the conditions
of
the last part of Theorem 2.1 are satisfied here, or to look
up the direct argument which is in the first version of this
paper
(\href{http://www.arxiv.org/abs/0812.3084}{arXiv:0812.3084}).
Thus (\ref{varVlim}) holds in its entirety.
\end{pf*}

The computations for $S$ are somewhat similar.
With $X_i$ denoting the indicator of the event that
$U_i$ is isolated,
\begin{eqnarray*}
\operatorname{Var}(S)
& = &
n \operatorname{Var}(X_1) + n(n-1) \operatorname{Cov}(X_1,X_2)
\\
& = & n (1 - \phi/n)^{n-1} \bigl(1 - (1 - \phi/n)^{n-1}\bigr)
+ n (n-1) \operatorname{Cov}(X_1,X_2).
\end{eqnarray*}
Since
$\operatorname{Cov}(X_1,X_2) = \mathbb E[X_1X_2] - \mathbb E[X_1]^2$,
provided $4\rho< n^{1/d}$ we can write
%
%e7.3 ###
%
\begin{eqnarray}\label{0521a}\quad
\operatorname{Var}(S)
& = & n (1 - \phi/n)^{n-1} \bigl(1 - (1 - \phi/n)^{n-1}\bigr)
\nonumber\\
& &{} + (n-1) \int_{B_{2 \rho} (\mathbf{0}) \setminus B_\rho(\mathbf{0})}
\biggl( 1- \frac{\rho^d \omega_d(|y|/\rho)}{n} \biggr)^{n-2} \, dy
\\
& &{} + n(n-1) \biggl( \biggl(1 - \frac{2^d \phi}{n} \biggr) \biggl(1 -
\frac{2\phi}{n} \biggr)^{n-2} - \biggl(1 - \frac{\phi}{n}
\biggr)^{2n-2}\biggr).\nonumber
\end{eqnarray}
\begin{pf*}{Proof of (\protect\ref{varWlim})}
For asymptotics as $n \to\infty$ with $\rho$ fixed,
by again using the MacLaurin expansion of $\log(1-x)$ we obtain
\begin{eqnarray*}
\biggl(1-\frac{2\phi}{n} \biggr)^{n-2} & = &
\exp\biggl((n-2) \biggl(-\frac{2\phi}{n}-\frac{2\phi^2}{n^2} +
O(n^{-3})\biggr) \biggr)
\\
& = &\exp\biggl(-2\phi+ \frac{ 4 \phi- 2 \phi^2}{n} + O(n^{-2})\biggr)
\end{eqnarray*}
and
\begin{eqnarray*}
\biggl(1-\frac{\phi}{n}
\biggr)^{2n-2}
& = &\exp\biggl((2n-2)\biggl(-\frac{\phi}{n}-\frac{\phi^2}{2n^2} +
O(n^{-3})\biggr) \biggr)
\\
& = &
\exp\biggl(-2\phi+ \frac{ 2 \phi- \phi^2}{n} +
O(n^{-2}) \biggr)
\end{eqnarray*}
and hence the last term in the right-hand side of (\ref{0521a})
is equal to
\begin{eqnarray*}
&&
n(n-1)
\exp(-2\phi)
\\
&&\quad{} \times\biggl( \biggl(1 - \frac{2^d \phi}{n}
\biggr)\exp\biggl(\frac{4 \phi- 2\phi^2}{n}\biggr)
- \exp\biggl(\frac{2\phi- \phi^2}{n} \biggr)
+ O(n^{-2}) \biggr)\\
&&\qquad= n(n-1) \exp(-2\phi) \biggl(- \frac{2^d \phi}{n} + \frac
{2 \phi}{n} -
\frac{\phi^2}{n} + O(n^{-2}) \biggr),
\end{eqnarray*}
so that
\begin{eqnarray*}
&& \lim_{n \to\infty} n^{-1}
\operatorname{Var}(S)
\\
&&\qquad = e^{-\phi} (1 - e^{-\phi})
- e^{-2 \phi}\bigl( (2^d -2) \phi+ \phi^2\bigr)
%& &
+ \int_{B_{2\rho}(\mathbf{0}) \setminus B_\rho(\mathbf{0})}
e^{- \rho^d \omega_d(|y|/\rho)} \,dy
\\
&&\qquad = e^{-\phi} - \bigl(1 + (2^d-2) \phi+ \phi^2\bigr)
e^{-2 \phi} + \rho^d \int_{B_{2}(\mathbf{0}) \setminus B_1(\mathbf{0})}
e^{- \rho^d \omega_d(|u|)} \,du,
\end{eqnarray*}
and since this limit is equal to $g_S(\rho)$ as
defined by (\ref{gWdef}),
we have proved the first
part of (\ref{varWlim}), namely, convergence to $g_S(\rho)$.

To complete the proof of (\ref{varWlim}), we need to
show that $g_S(\rho) > 0$. This can be done by
the same arguments as for the
% as was used in the
proof of (\ref{varVlim}).
% already given.
Hence,
(\ref{varWlim}) holds in its entirety.
\end{pf*}
\begin{pf*}{Proof of Theorem \protect\ref{thmlimsup}}
It remains only to prove
(\ref{main2}),
(\ref{main1})
and (\ref{SLB}).
%It is immediate from the
By definition
% \eq{etaVinf}, \eq{etaWinf}, \eq{etaVdef} and \eq{etaYdef}
%that
$ \eta_V(\rho)
= \lim_{n \to\infty} \eta_V(n,\rho) $
and
$ \eta_S(\rho)
= \lim_{n \to\infty} \eta_S(n,\rho)$.
Then
(\ref{main2})
follows at once from Theorem \ref{thm2}, along with
(\ref{meanlim}) and (\ref{varVlim}).
% and \eq{etaVinf}.
Similarly,
(\ref{main1})
follows at once from Theorem \ref{thm1} along with (\ref{meanlim}),
and
(\ref{varWlim}).
% and \eq{etaWinf}.
% $\qed$

Finally, we demonstrate the asymptotic lower bound (\ref{SLB}).
For any random variable $X$,
let $F_X$ denote its cumulative distribution function
and
let
$f_X$ denote its probability density function (if it has one).
Let $\varepsilon\in(0,1)$.
Set
\[
t_1 : = \frac{ [\mu_S] -\mu_S }{\sigma_S} ;\qquad
t_2 := \frac{ [\mu_S] -\mu_S +1 - \varepsilon}{\sigma_S} .
\]
Here $[\cdot]$ denotes integer part, so that $|t_i| \leq\sigma_S^{-1}
$ for $i=1,2$.
By the unimodality of the standard normal density,
%
%e7.4 ###
%
\begin{eqnarray}\label{1031a}
F_Z(t_2 ) - F_Z(t_1) & \geq&
(t_2-t_1) \min(f_Z(t_1),f_Z(t_2)) \nonumber\\[-8pt]\\[-8pt]
& \geq&
(1 - \varepsilon) \sigma_S^{-1} f_Z(\sigma_S^{-1}).\nonumber
\end{eqnarray}
On the other hand, since $S$ is integer-valued, $F_{(S-\mu_S)/\sigma_S}(t_1) $
is equal to\break
$F_{(S-\mu_S)/\sigma_S}(t_2)$,
so that by (\ref{1031a})
\[
%d_K (\LL((S-\mu_S)/\sigma_S), \LL(Z) )
D_S
\geq(1/2)
(1 - \varepsilon) \sigma_S^{-1} f_Z(\sigma_S^{-1}).
\]
Scaling by $n^{1/2}$,
% then
letting $n \to\infty$, using (\ref{varWlim}) and
%then
letting $\varepsilon\to0$
yields (\ref{SLB}).
\end{pf*}

To conclude,
we compute some numerical values for
the asymptotic upper bounds
% $\delta_W(\rho)$ and $\delta_V(\rho)$
appearing in
(\ref{main2}) and (\ref{main1}).
% Theorem \ref{thmlimsup}.
For this we need
to compute $J_{r,d}(\rho)$ defined by (\ref{Jdef}) (for $r=1$ and $r=2$),
and for this in turn,
we need to compute $\omega_d(u) $, the volume of the
%intersection
union of two unit balls in $d$-space whose centers are at points
($x,x'$ say) distance $u$ apart ($u \leq2 $).
%Denote these balls by $B$ and $B'$ respectively.
Clearly,
$ \omega_1(u) = 2 + u$, and
generalizing (6) of \cite{Moran1}
to arbitrary $d \geq2$, we have
%
%e7.5 ###
%
\begin{equation}\label{omMoran}
\omega_d(u) = \pi_d + \pi_{d-1} \int_0^u \bigl( 1- (t/2)^2 \bigr)^{(d-1)/2}
\,dt,\qquad
d \geq2.
\end{equation}
%
%In particular,
% \omega_3(u) &= & \pi( (4/3) + r - (1/12)r^3 ).
%By \eq{Jdef}
%so that by \eq{Jdef}
% one may therefore derive
%J_{r,1}(\rho) & = & 2 e^{-2 \rho} \int_0^r \exp(-\rho t) dt =
% 2 e^{-2 \rho} (1 - e^{-\rho r} ) /\rho; \\
%J_{r,2}(\rho)
%& = &
% 8 \pi e^{- \pi\rho^2}\int_0^{r/2} \exp(- 2 \rho^2
%[\arcsin(t) + t (1 - t^2)^{1/2} ] ) t dt ;
%J_{r,3}(\rho) & = &
%4 \pi e^{-(4/3) \pi\rho^3 }
% \int_0^r \exp( \pi\rho^3 (\frac{s^3}{12} -
%s ) ) s^2 ds.
Using the preceding formulae,
%in \eq{deltaV} and \eq{deltaS},
we have computed numerical values for the asymptotic
upper bounds in Theorem \ref{thmlimsup},
for the cases with $\rho=1$ and $d \leq3$. These are as follows to
five significant figures, where $\delta_V(\rho)$ denotes the right-hand
side of (\ref{main2}) and $\delta_S(\rho)$ denotes the right-hand
side of (\ref{main1}):
% {\bf values need updating!}
%
%$$
% \begin{array}{lll}
%1.1391 \times10^5 & \mathrm{if} & d=1 \\
%1.8193 \times10^7 & \mathrm{if} & d=2 \\
%3.2329 \times10^9 & \mathrm{if} & d=3
% \begin{array}{lll}
%2.7270 \times10^4 & \mathrm{if} & d=1 \\
%6.3096 \times10^5 & \mathrm{if} & d=2 \\
%1.4482 \times10^7 & \mathrm{if} & d=3
%$$
%
%Optimized bounds
%$$
% \begin{array}{lll}
%6.5785 \times10^3 & \mathrm{if} & d=1 \\
%8.6418 \times10^5 & \mathrm{if} & d=2 \\
%1.4455 \times10^8 & \mathrm{if} & d=3
% \begin{array}{lll}
%2.4132 \times10^3 & \mathrm{if} & d=1 \\
%5.2628 \times10^4 & \mathrm{if} & d=2 \\
%1.1789 \times10^6 & \mathrm{if} & d=3
%$$
%Doubly Optimized bounds
%
\begin{eqnarray*}
\delta_V (1) &=& \cases{
6.4252 \times10^3, &\quad if $d=1$, \cr
8.6212 \times10^5, &\quad if $d=2$, \cr
1.4451 \times10^8, &\quad if $d=3$,}
\\
\delta_S (1) &=& \cases{
2.1024 \times10^3, &\quad if $d=1$, \cr
4.6833 \times10^4, &\quad if $d=2$, \cr
1.0578 \times10^6, &\quad if $d=3$.}
\end{eqnarray*}

\section*{Acknowledgments}
This work was partly done at meetings in
% taking place at
the Institute for Mathematical Sciences at the National University
of Singapore, whose
% hospitality and
support we gratefully acknowledge.
We also thank Joseph Yukich for conversations which stimulated
our initial interest in this topic.
We are indebted to the referee for the idea for Lemma \ref{LGthm},
and other helpful suggestions.
%in 2008 and 2009, at
% It is a pleasure to acknowledge the value of conversations
%with
%We began
%this work while visiting the Institute for Mathematical Sciences
%at the National University of Singapore, whom
% we thank for their hospitality.

% imsref loaded by lrinkeviciute, 2009-10-07 14:52:01
%

%
\printaddresses

\end{document}